\documentclass[red,11pt,a4paper]{article}

\usepackage{cite}

\usepackage{amsmath}
\usepackage{amscd}
\usepackage{amssymb}
\usepackage{latexsym}
\usepackage{color}
\numberwithin{equation}{section}
\usepackage[normalem]{ulem}

\renewcommand{\div}{{\rm div}\,}

  \newcommand{\vrp}{{\vr^+}}
    \newcommand{\vrm}{{\vr^-}}    
\def\calK{{\mathcal K}}

\def\d{\partial}

\newcommand{\T}{\mathbb{T}}

\renewcommand{\d}{\partial}

\newcommand{\lap}{\Delta}

\newcommand{\Div}{\operatorname{div}}

\newcommand{\pf}{{\noindent\it Proof.~}}
\newcommand{\Ov}[1]{\overline{#1}}

\newcommand{\vr}{\varrho}

\newcommand{\vu}{\vc{u}}

\newcommand{\vc}[1]{{\bf #1}}

\newcommand{\Grad}{\nabla}

\newcommand{\pt}{\partial_{t}}
\newcommand{\ptb}[1]{\partial_{t}(#1)}
\newcommand{\Dt}{\frac{ d}{dt}}

\newcommand{\dx}{{\rm d} {x}}
\newcommand{\dy}{{\rm d} {y}}
\newcommand{\dt}{{\rm d} t }

\newcommand{\dxdt}{\dx \,\dt}
\newcommand{\lr}[1]{\left( #1 \right)}
\newcommand{\intO}[1]{\int_{\T^d} #1 \ \dx}
\newcommand{\iintO}[1]{\int_{\T^{2d}} #1 \ \dx\,\dy}

\newcommand{\intOB}[1]{\int_{\T^d} \left( #1 \right) \ \dx}

\newcommand{\intT}[1]{\int_0^T #1 \ \dt}
\newcommand{\intTO}[1]{\int_0^T\!\!\!\! \int_{\T^d} #1 \ \dxdt}

\newcommand{\eq}[1]{\begin{equation}
\begin{split}
#1
\end{split}
\end{equation}}
\newcommand{\eqh}[1]{\begin{equation*}
\begin{split}
#1
\end{split}
\end{equation*}}
\newcommand{\ep}{\varepsilon}

\newcommand{\R}{\mathbb{R}}
\newtheorem{thm}{Theorem}[section]
\newtheorem{lemma}[thm]{Lemma}
\newtheorem{prop}[thm]{Proposition}

\newtheorem{rmk}{Remark}
\newtheorem{cor}[thm]{Corollary}

\title{Finite-Energy Solutions \\  for  Compressible Two-Fluid Stokes System }

\author{Didier Bresch\thanks{LAMA CNRS UMR 5127, University of Savoie Mont-Blanc, Bat. Le Chablais, Campus scientifique,73376 Le Bourget du Lac, France, E-mail: \texttt{didier.bresch@univ-smb.fr}}, Piotr B. Mucha\thanks{Institute of Applied Mathematics and Mechanics, University of Warsaw,  ul. Banacha 2, 02-097 Warszawa, Poland, E-mail: \texttt{p.mucha@mimuw.edu.pl}}, Ewelina Zatorska\thanks{University College London, Department of Mathematics,  Gower Street, London WC1E 6BT, United Kingdom. Email: \texttt{e.zatorska@ucl.ac.uk}}}

\begin{document}
 \maketitle 

\abstract
We prove existence of global in time weak solutions to a compressible two-fluid Stokes system  with a single velocity field and algebraic closure for the pressure law. The constitutive relation involves  densities of both fluids through  an implicit function. The system appears to be outside the class of problems that can be treated using the classical Lions-Feireisl approach. Adapting the novel compactness tool developed by the first author and P.--E. {\sc Jabin} in the mono-fluid compressible Navier-Stokes setting, we first prove the weak sequential stability of solutions. Next, we construct weak solutions via unconventional approximation using  the Lagrangian formulation, truncations and stability  result of trajectories for rough velocity fields.

\medskip


\section{Introduction}
Multi-component fluid models arise in various applications including studies of water wave impact on  coastal structures (violent aerated flows) \cite{DiDuGh}, petroleum industry \cite{Evje08}--\cite{Evje16}, cancer cell migration in compressible media \cite{Ev}, or turbulent mixing in nuclear industry, reactive flows, propulsion and sprays \cite{De}, to name only few. Classical derivation of multi-fluid models begins writing the equilibrium equations for each component of the flow at the microscopic level  with small scale interfaces. 
The second step is to  perform a volume averaging under suitable  closure assumptions. Averaged models  bypass the local geometrical complexity of the interphase at the cost of including new variables --
the volumetric rate of presence of each fluid/phase, see \cite{De}. This formal derivation can be found in the monographs of {\sc M. Ishii} and {\sc T. Hibiki} \cite{IsHi},  and of {\sc D. Drew} and {\sc S. L. Passman} \cite{DrPa}. 
Mathematically rigorous derivation of several models from mono-fluid systems  may be found in \cite{BrHi}, \cite{BrHi2}, \cite{BrHu}, \cite{AmZl}, \cite{Se}.  The  reader is also referred to the recent chapter \cite{BrDeGhGrHi} for discussion on modelling and mathematical studies of multi-fluid systems in the compressible setting.

\smallskip 

    In the present paper, we analyze a bi-fluid compressible system  in the semi-stationary Stokes regime. 
    We assume a common velocity field and pressure for both fluids (algebraic pressure closure). Our system of equations reads:
\begin{subequations}\label{spec}
\begin{align}
&\ptb{\alpha^\pm\vr^\pm}+\Div(\alpha^{\pm}\vr^\pm\vu)=0,\label{spec_1}\\
&-\mu\lap\vu - (\lambda+\mu) \nabla \Div \vu +\Grad p=\vc{0},\\
&\alpha^++\alpha^-=1,\\
& p = p^+ = p^- ,
\end{align}
\end{subequations}
with constant shear and bulk viscosities $\mu$ and $\lambda$ such that  $\lambda+2\mu>0$ {and $\mu>0$}. The unknowns of the system \eqref{spec} are  the volumetric rates of presence of fluid $+$ {and} $-$,  $\alpha^+,\alpha^-$, respectively,  with  
\eq{\label{alpha}
0\le \alpha^\pm \le 1,
}
the two mass densities $\vrp$, $\vrm$, and the common velocity field $\vu$. By $p^+$, $p^-$ we denote the internal barotropic pressures for each fluid with the explicit form:
\eq{\label{CP}
p^+=a^+{\vr^+}^{\gamma_+}, \quad p^-=\quad a^-{\vr^-}^{\gamma-}, 
}
where $a^\pm>0$, $\gamma^\pm>1$ are given constants. 
The purpose of this paper is to prove the existence of solutions "\`a la {\sc Leray}" (finite energy) to this system on the $d$-dimensional torus $\T^d$, $d=2,3$,  under the following constraint 
\eq{\label{av}
\intO{ \vu(t,x)}=0 \quad \mbox{ for } \  t\in (0,T),}
and with the initial conditions
\eq{\label{Ic_1}
\alpha^+ \vr^+ \vert_{t=0} = R_0, \quad \alpha^- \vr^- \vert_{t=0} = Q_0, \qquad 
R_0 \ge 0, \quad Q_0 \ge 0.\\
}
Moreover, we ask for the following compatibility condition for the initial data:
$$\alpha^+ \vert_{t=0} = \alpha^+_0,\qquad \alpha^- \vert_{t=0} = \alpha^-_0$$
\eq{\label{Ic_2}
\alpha^+_0 + \alpha^-_0 = 1, \qquad \alpha^\pm_0 \ge 0, \qquad
p^+(\vr^+_0) = p^-(\vr^-_0).
}
\begin{rmk}\label{Rmk_1}
If $\vr_0^+$ is nonzero, then $\alpha_0^+ = R_0/\vr_0^+$. Because of the algebraic pressure law closure which provides $\vr_0^-= (a_+/a_-)^{1/\gamma_-} (\vr_0^+)^{{\frac{\gamma^+}{\gamma^-}}} \not = 0$ we can also define  $\alpha_0^- = Q_0/\vr_0^-$.
If, on the other hand, $\vr_0^+ = \vr_0^- = 0$, we choose $\alpha_0^\pm  = 1/2$ for each phase.  
\end{rmk}
 
 Our main result for system \eqref{spec} reads:
 \begin{thm}. \label{Mainbifluid}
Let  $\gamma^\pm>1$, {$\gamma^+\neq \gamma^-$,} $a^\pm>0$, $\lambda + 2\mu>0$, {$\mu>0$}, and let the initial data  \eqref{Ic_1} satisfy \eqref{Ic_2}. 
Assume that 
$$\intO{ (R_0^{\gamma^+} + Q_0^{\gamma^-}) }< \infty, \qquad
   0< \intO{ R_0 } <  \infty, \qquad
   0 < \intO{ Q_0 } < \infty.$$ 
Then there exists a global weak solution $(\alpha^\pm, \vr^\pm,\vu)$ of system \eqref{spec}--\eqref{av} 
satisfying
$$\alpha^\pm\vr^\pm \in  L^\infty(0,T;L^{\gamma^\pm}(\T^d)) \cap L^{2\gamma^\pm}((0,T)\times \T^d) \cap {\cal C}([0,T]; L^{\gamma^\pm}(\T^d)),
$$
$$\vr^\pm \in  L^\infty(0,T;L^{\gamma^\pm}(\T^d)) \cap L^{2\gamma^\pm}((0,T)\times \T^d) \cap {\cal C}([0,T]; L^{\gamma^\pm}(\T^d)),
$$
$$ \vu \in L^2(0,T; H^1(\T^d)),$$
where Equations \eqref{spec} are satisfied in ${\cal D}'((0,T)\times \T^d)$ and the initial conditions \eqref{Ic_1}
are satisfied in ${\cal D}'(\T^d)$. Moreover, the equations \eqref{spec_1} are satisfied in the renormalized sense.
\end{thm}

   Note  that if $a^+=a^-=1$ and $\gamma^+ = \gamma^-$, system \eqref{S} is reduced to the 
 semi-stationary Stokes version of  \eqref{entropytransport2}. In this case our technique provides stronger results that the classical technique of Lions. Indeed comparing \cite[Theorem 2]{MaMiMuNoPoZa} (see also \cite[Theorem 3.1]{Michalek}) with the proof of Theorem \eqref{Thm:main} we see that we are able to prove compactness of sequences approximating $R,Q$ and $Z$, while the results from \cite{MaMiMuNoPoZa} provide strong convergence of sequence approximating $Z$, but only a weak convergence of $\vr_n\to\vr$.
 To get Theorem \ref{Mainbifluid}, we first prove global existence of weak solutions for a reformulation of the system \eqref{spec}.  Introducing the notation
$$R =\vr^+\alpha^+,\qquad Q =\vr^-\alpha^-, \qquad Z= \vr^+,$$
we check that the pressure $p$ is expressed in terms of $R,Q$. In fact we have
\begin{equation}\label{pZ}
 p=P(R,Q)= a^+ {Z}^{\gamma_+},
\end{equation}
for $Z=Z(R,Q)$ such that 
\eq{\label{TZ}
\lr{\frac{a^-}{a^+}}^{1/\gamma^-} Q = \lr{1-\frac{R}{Z} } Z^\gamma,\quad  \mbox{with}\quad \gamma=\frac{\gamma_+}{\gamma_-},
}
and
\eq{\label{RleqZ}R\leq Z.}
The system \eqref{spec}-\eqref{av}  can be therefore transformed to the following form
\begin{subequations}\label{S}
\begin{align}
  &\pt R+ \Div(R\vu )=0,\label{SR}\\ 
  &\pt Q + \Div(Q\vu )=0 ,\label{ST}\\
&- (\lambda + 2 \mu)  \Div \vu + a^+\lr{Z(R,Q)^{\gamma^+}  - \{ Z(R,Q)^{\gamma^+} \} }=0, \label{Div}\\
&{\rm rot} \, \vu = 0, \quad \intO{ \vu(t,x) }= 0, \label{rotu}
\end{align}
\end{subequations} 
where $\{f \}= \bigl( \int_{\T^d} f(x) \, dx\bigr)/|\T^d|$ and $Z$ is related to $R$ and $Q$ through the non-explicit formula \eqref{TZ}.
The initial conditions for the previous system reads
\eq{\label{IcNew}
R \vert_{t=0} = R_0, \quad Q \vert_{t=0} = Q_0,\qquad  R_0 \ge 0, \quad Q_0 \ge 0,
}
with the following compatibility condition on  $ Z\vert_{t=0} = Z_0$
\eq{\label{RelationInit}
\lr{\frac{a^-}{a^+}}^{1/\gamma^-} Q_0 = \lr{1-\frac{R_0}{Z_0} } Z_0^\gamma, 
\qquad  \hbox{ with } \qquad R_0 \le Z_0.
}
For system \eqref{S} we prove the following:
\begin{thm}\label{Thm:main}
Let $\gamma^\pm>1$, {$\gamma^+\neq \gamma^-$}, $a^\pm>0$, $\lambda + 2\mu>0$, {$\mu>0$}, and let the initial data be given by \eqref{IcNew} with $Z_0$
defined through \eqref{RelationInit}. Further, assume that
\eq{\label{eq_ini}
\intO{ (R_0^{\gamma^+} + Q_0^{\gamma^-}) }< \infty,  \qquad
   0<  \intO {R_0} < \infty, \qquad  0< \intO {Q_0} < \infty.
}
Then there exists $(R,Q,Z,\vu)$ -- a global in time weak solution to the system \eqref{TZ}--\eqref{S}  for $(t,x)\in (0,T)\times \T^d$,
 such that
$$R  \in L^\infty(0,T;L^{\gamma^+}(\T^d)) \cap L^{2\gamma^+}((0,T)\times \T^d) \cap {\cal C}([0,T]; L^{\gamma^+}(\T^d)),
$$
$$
   Q \in L^\infty(0,T;L^{\gamma^-}(\T^d)) \cap L^{2\gamma^-}((0,T)\times \T^d) \cap {\cal C}([0,T]; L^{\gamma^-}(\T^d))
$$
$$Z  \in L^\infty(0,T;L^{\gamma^+}(\T^d)) \cap L^{2\gamma^+}((0,T)\times \T^d) \cap {\cal C}([0,T]; L^{\gamma^+}(\T^d)),
$$
$$\vu \in L^2(0,T;H^1(\T^d)),
$$
where Equations \eqref{pZ}-\eqref{rotu} are satisfied in ${{\cal D}'((0,T)\times \T^d})$, and the initial conditions \eqref{IcNew}  with 
the constraint \eqref{RelationInit} are satisfied in  ${\cal D}'(\T^d)$. Moreover, equations \eqref{SR} and \eqref{ST} are satisfied in the renormalized sense.
\end{thm} 
 
Our paper provides the first proof of existence of global-in-time weak solutions for a bi-fluid system with constant viscosities under algebraic pressure law closure in
physical dimensions 2 and 3. The only other available results treat the density dependent viscosities with two velocity fields, see \cite{BrDeGhGr}, and \cite{BrHuLi} (see also \cite{Evje16}--\cite{Evje16a} for specific
linear pressure laws). The main difficulty in analysis of the system 
\eqref{spec} is due to, roughly speaking, complex form of the pressure. Indeed, by $\alpha^++\alpha^-=1$
 and \eqref{CP}, the pressure $p$ depends on  $\alpha^+\vr^+$ and $\alpha^-\vr^-$ ($R$ and $Q$, respectively) in  nonlinear implicit, see the relation \eqref{TZ}. It causes that the nowadays standard approach developed by {\sc P.-L. Lions} \cite{Lions2} and {\sc E. Feireisl} \cite{EF2001}, see also \cite{MaMiMuNoPoZa},  seems to be not applicable in all generality. 
 Therefore, we adapt a brand new technique  from \cite{BJ,BJ_short, BrJa} that could be used
 to cover more complicated case in a future work  for instance with viscosities depending
on the volumic rates $\alpha^\pm$.   
  By reformulating system \eqref{spec} in terms of the pressure argument $Z$ and one of the conserved  quantities, for example $R=\alpha^+\vr^+$, we show that the new technique provides compactness for sequences approximating both unknowns. 

Let us briefly discuss various contributions related to global weak solutions to the study of compressible fluid equations with intricate pressure law.
\medskip

\noindent {\bf Mono-fluid systems.} Note that putting $\alpha^+ =1$ in \eqref{spec}, we get the usual semi-stationary compressible Stokes system
\eq{
&\pt \vr +\Div(\vr\vu)=0,\\
&-\mu\lap\vu - (\lambda+\mu) \nabla \Div \vu +\Grad p(\vr)=\vc{0},
\label{Quasi-Stokes}}
that has been studied for instance in \cite{Lions2} with a monotone pressure law $p(\vr) = a \vr^\gamma$ with $\gamma>1$, and
more recently in \cite{BJ_short} with non-monotone, locally Lipschitz, pressure law $p(\vr)$, such that $p(0)=0$
and 
$$C^{-1} \vr^\gamma - C \le p(\vr) \le C\vr^\gamma + C, \qquad
   |p'(s)| \le \overline p s^{\gamma-1},
$$
for some constants $C>0$, $\overline p >0$ and $\gamma>1$.
\medskip

\noindent {\bf Compressible systems with two continuity equations.}  Two-components compressible systems have been studied in the density dependent viscous case in \cite{BrDeGhGr}, \cite{BrHuLi} for instance. 
Concerning the constant viscosity case, 
the existence of weak solutions to the two--phase model 
{\eq{
&\pt \vrp + \Div(\vrp\vu)=0\\
&\pt\vrm + \Div(\vrm \vu) = 0\\
& \pt ([\vrp+\vrm]\vu)+ \Div ([\vrp+\vrm] \vu\otimes \vu) + \nabla p(\vrp,\vrm) - \mu \Delta \vu - (\lambda + \mu) \nabla\Div \vu = \vc{0}
\label{two-phaseAV}}}
was recently proven by A. {\sc Vasseur}  {\it et al.} in \cite{VaWeYu}, for the pressure law equal to {$p(\vrp,\vrm) = \lr{\vrp}^\gamma + \lr{\vrm}^\alpha$} with $\gamma>9/5$ and $\alpha \ge 1$. For the existence of strong solutions close 
to equilibrium we refer, for example, to \cite{Hao12}. For the extensive analysis of two-component models in the one-dimensional setting we refer to papers of {\sc S. Evje}  {\it et al.} \cite{Evje08, Evje15} in the constant viscosity case
and with specific linear pressure law.

In another recent paper of D. {\sc Maltese}  {\it et al.} \cite{MaMiMuNoPoZa}, see also \cite{Michalek, FeKlNoZa}, the authors considered the system used in the geophysical flow modelling
\eq{
&\pt \vr + \Div(\vr \vu)=0\\
&\pt s  + \vu\cdot\Grad s  = 0\\
& \pt (\vr \vu) + \Div (\vr  \vu\otimes \vu)+ \nabla p(\vr,s) - \mu \Delta \vu - (\lambda + \mu) \nabla\Div \vu = \vc{0}
\label{entropytransport},}
where $s$ denotes the entropic variable, with the pressure law given by $p(\vr,s) = \vr^\gamma {\cal T}(s)$ with $\gamma>1$, $s>0$,
and ${\cal T}(\cdot)$ a given smooth and strictly monotone function. They proved the existence of weak solution to the following reformulation of \eqref{entropytransport}  
\eq{\label{entropytransport2}
&\pt \vr + \Div(\vr \vu)=0\\
&\pt Z  + \Div(Z \vu) = 0\\
& \pt (\vr \vu) + \Div (\vr  \vu\otimes \vu)+ \nabla Z^\gamma - \mu \Delta \vu - (\lambda + \mu) \nabla\Div \vu = \vc{0},}
where $Z$ denotes $\vr[{\cal T}(s)]^{1/\gamma}$, after which they proved the equivalence between solutions to systems \eqref{entropytransport} and \eqref{entropytransport2} for  $\gamma\geq \frac95$. For analysis and numerical simulations of system \eqref{entropytransport2} with the so-called congestion constraint we refer to \cite{DeMiZa, DeMiNaZa}.
 
 \medskip
Note that the two systems mentioned before \eqref{two-phaseAV} and \eqref{entropytransport} include pressure laws which are monotone with respect to variables satisfying continuity equations. This allows the authors to adapt the tools already developed by {\sc P.--L. Lions} \cite{Lions1, Lions2} and {\sc E.~Feireisl} \cite{EF2001} for mono-fluid systems.  It is not the case for the two-fluid Stokes system \eqref{S}. For this system, although $\partial_RZ,\ \partial_Q Z\geq 0$, we do not have the crucial property
$$\Ov{Z(R,Q)}\geq Z(R,Q),$$
for $\Ov{Z(R,Q)}$ denoting the weak limit of  $Z_n(R_n,Q_n)$ when the approximation parameter $n$ goes to $+\infty$.
Looking at the reformulation \eqref{S}, the pressure argument $Z$ is given implicitly in terms of $R$ and $Q$. As a consequence, $Z$ does not satisfy a simple continuity equation, {but} we have
\eqh{
\pt Z+\Div(Z\vu)+\frac{(1-\gamma) (Z-R)Z}{ \gamma(Z-R)+R}\Div\vu=0.}
The additional friction term $\frac{(1-\gamma) (Z-R)Z}{ \gamma(Z-R)+R}\Div\vu$ causes that 
its seems that we have no compensation of compactness between $\Div\vu$ and the pressure $p$. Continuity equation with production term in the Navier-Stokes type of system has been recently investigated by {\sc N. Vauchelet} and {\sc E. Zatorska} \cite{VauchZator}.
In this case, the so-called effective viscous flux equality does not imply a strong convergence of the sequence approximating $Z$. 

In this work we  show that the recent development proposed by {\sc D. Bresch} and {\sc P.--E. Jabin} \cite{BJ} may be adapted to treat the bi-fluid system \eqref{S}. Our  work therefore provides a  generalization  of this result to the pressure law that depends on two transported quantities, as in semi-stationary compressible Stokes system. 

In \cite{BJ}, the authors explain how to handle the non-monotone truncated pressure in the heat-conducting Navier-Stokes system. In this system the pressure  depends on two variables: the density and the temperature. The density satisfies the continuity equation, while the temperature is given by the heat equation, and hence some properties providing  compactness in space  of the second unknown are available.  Our result {in this paper} covers the pressure laws depending on two quantities without knowing any {\it a-priori} compactness in space for any of them.

\medskip

In our proof we rescale the unknowns  and the  viscosity  coefficient $\lambda + 2 \mu$ so that  $a^+=a^-=1$. Since we keep $\gamma^+\neq \gamma^-$ this assumption does not lead  to loss of generality. For the sake of brevity, we will always consider $\gamma^+\le \gamma^-$, equivalently $\gamma\leq 1$. However,  due to the symmetry of the problem, the result will remain in force also if $\gamma^+>\gamma^-$. 

\medskip

The paper is divided into two parts:

\smallskip
\noindent {\bf Part I.} In Section \ref{Sect:2}, we first prove energy estimates and extra integrability properties on the solutions of the system \eqref{S}. We also study the nonlinear relation
between $Q$, $R$ and $Z$ and present the equation satisfied by $Z$. Then, in Section \ref{Sect:3},  we prove the weak sequential stability of solutions to \eqref{S}. This means that the hypothetical sequence of sufficiently smooth solutions
$\{R_n,Q_n,Z_n,\vu_n\}_{n=1}^\infty$ satisfying  the energy and extra integrability estimates uniformly w.r.t. $n$, has a limit when $n\to\infty$,
that is a weak solution to \eqref{S}.

\smallskip
   
\noindent {\bf Part II.} In Section \ref{Sect:4}, we construct the approximate solutions and show that they converge to solutions of system \eqref{S}.  The starting point is the Lagrangian reformulation of \eqref{S} with truncation of the pressure. Using  {\sc G. Crippa} and {\sc C. De Lellis}' stability of the flow  result (see \cite{CrDe} and \cite{CoCrSp}) we show that approximate solutions constructed in
the Lagrangian coordinates define suitable approximate solutions of the system in the Eulerian coordinates. These solutions satisfy the uniform bounds requested  in the first part of the paper. It is worth to emphasize that our construction does not introduce any parabolic regularization of the continuity equation commonly used in compressible setting. In a sense it is similar to construction of regular solutions \cite{Mucha2003, MuchaZaj2002, MuchaZaj2006}.

\section{Preliminary observations}\label{Sect:2}
Here we provide basic a-priori estimates for the sequence of solutions $\{R_n, Q_n, Z_n, \vu_n\}_{n=1}^{\infty}$, uniformly with respect to $n$. We assume that for any $n\geq 1$ $(R_n, Q_n, Z_n, \vu_n)$ is a
smooth solution to \eqref{S}, defined on $(0,T)\times\T^d$. We drop the index $n$, when no confusion can arise, and we recall that we assume $\gamma\leq 1$.
Moreover, the results from this section do not depend on the value of the viscosity coefficient if only $\lambda+2\mu>0$. Therefore, without loss of generality we take
{$$\lambda+2\mu=1.$$}

\begin{lemma}\label{Lemma_posit}
 Let $R,Q,\vu$ be sufficiently smooth solutions to (\ref{S}), then 
 \begin{equation}\label{RQ0}
 0\leq R,Q \quad \mbox{and}\quad R,Q\in L^\infty(0,T; L^1(\T^d)).
\end{equation}
Moreover, assuming $R$, $Q \ge 0$ given,  there exists a unique $Z$ solving \eqref{TZ} and \eqref{RleqZ}.
\end{lemma}
\pf Integrating  equation \eqref{SR} over $\T^d$ we deduce that
\begin{equation*}
\Dt \intO{R}=0,
\end{equation*}
i.e.  $\intO{R_{0}(x)}=M_R$, implies $\intO{R(t,x)}=M_R$ for any $t\in[0,T]$.\\
Moreover, since $R$ is smooth and $R_0\geq0$, we have the following estimate
\begin{equation}\label{ch2:ron0}
R(\tau,x)\geq \inf_{x\in\T^d}R_0(x)\exp\left(-\int_{0}^\tau\|\Div\vu\|_{L^\infty(R)}{\rm d}t\right),
\end{equation}
in particular $R\geq0$. Repeating the same procedure for $Q$ we obtain \eqref{RQ0}.

For fixed nonnegative $R$ and $Q$, we find
a candidate $Z\ge R$ satisfying the equation \eqref{TZ}, i.e.
$$f_{R,Q}(Z) = Z^\gamma - R Z^{\gamma-1} - Q = 0.$$
Note that $f_{R,Q}(R) = - Q \le 0$ and  $\partial_Z f_{R,Q}(Z)= \gamma Z^{\gamma-2}(Z-R) + R Z^{\gamma-2} $. 
Therefore in the range $Z\geq R$ there exists a unique $Z$ solving \eqref{TZ}.  $\Box$

\bigskip

Our next goal is to derive estimates for $R,\,Q,\, Z$, and $\vu$ following from boundedness of the energy associated with system \eqref{S}.
\begin{lemma} \label{Lemma_est}
 Let $R,Q,\vu$ be sufficiently smooth solutions to (\ref{S}), then the following estimates are valid
 \begin{equation}\label{est_ener}
  \sup_{t<T} \intO{ (Z^{\gamma^+}+R^{\gamma^+} + Q^{\gamma^-}) } +  \intTO{|\Grad\vu|^2 }\leq C,
 \end{equation}
and 
\begin{equation}\label{est_bog}
 \intTO{ (Z^{2\gamma^+}+R^{2\gamma^+} + Q^{2\gamma^-})} \leq C(1+T),
\end{equation}
for any $T>0$.
\end{lemma}

\pf 
First note that from $\intO{Z_0^{\gamma^+}} <  \infty$ and \eqref{RelationInit} it follows that
\eq{\label{R0Q0}
\intO{R_0^{\gamma^+}} <  \infty,\quad \intO{Q_0^{\gamma^-}} <  \infty.}
On the opposite, assuming that \eqref{R0Q0} holds then using the H\"older inequality and \eqref{RelationInit}  it follows that $\intO{Z_0^{\gamma^+}} <  \infty$.
  We next define $\alpha$ as a solution to
  {\eqh{
  \alpha=\left\{
  \begin{array}{ll}
  \frac{R}{Z}&  \hbox{ if }  \  Z\not = 0\\
  \frac{1}{2} & \hbox{otherwise}.
  \end{array}
  \right.
  } }
 By \eqref{RleqZ} we deduce that $0\leq \alpha\leq 1$. 
From \eqref{TZ} we  have that 
{\eq{\label{equiv}
Z^{\gamma^+}=\lr{\frac{Q}{1-\alpha}}^{\gamma^-} \ \hbox{if}\ \alpha\neq1\quad \text{ or equivalently}\quad  Z^{\gamma^{+}}=\lr{\frac{R}{\alpha}}^{\gamma^+}\ \hbox{if}\ \alpha\neq 0.} }
Therefore, the gradient of the pressure can be written as
\eq{
\Grad P&=\Grad Z^{\gamma^+}=\alpha\Grad \lr{\frac{R}{\alpha}}^{\gamma^+}+(1-\alpha)\Grad\lr{\frac{Q}{1-\alpha}}^{\gamma^-}
\\
&=\frac{\gamma^+}{\gamma^+-1}R \Grad \lr{\frac{R}{\alpha}}^{\gamma^+-1}+\frac{\gamma^-}{\gamma^--1}Q\Grad \lr{\frac{Q}{1-\alpha}}^{\gamma^--1}.
}
{Multiplying  the last term on the l.h.s. of \eqref{Div} by $-\Div\vu$, and integrating by parts, we obtain}
\eqh{
&-\intO{(P-\{P\})\Div\vu}=\intO{\Grad P\cdot\vu}\\
&=\intO{\lr{\frac{\gamma^+}{\gamma^+-1}R \Grad \lr{\frac{R}{\alpha}}^{\gamma^+-1}+\frac{\gamma^-}{\gamma^--1}Q\Grad \lr{\frac{Q}{1-\alpha}}^{\gamma^--1}}\cdot\vu}\\
&=-\frac{\gamma^+}{\gamma^+-1}\intO{\Div(R\vu)\lr{\frac{R}{\alpha}}^{\gamma^+-1}}
-\frac{\gamma^-}{\gamma^--1}\intO{\Div(Q\vu)\lr{\frac{Q}{1-\alpha}}^{\gamma^--1}}\\
&=\frac{\gamma^+}{\gamma^+-1}\intO{\pt R\lr{\frac{R}{\alpha}}^{\gamma^+-1}}+\frac{\gamma^-}{\gamma^--1}\intO{\pt Q\lr{\frac{Q}{1-\alpha}}^{\gamma^--1}}\\
&=\frac{1}{\gamma^+-1}\Dt\intO{ \lr{\frac{R}{\alpha}}^{\gamma^+}\alpha}+\intO{ \lr{\frac{R}{\alpha}}^{\gamma^+}\pt \alpha}\\
&\qquad+\frac{1}{\gamma^--1}\Dt\intO{ \lr{\frac{Q}{1-\alpha}}^{\gamma^-}(1-\alpha)}-\intO{ \lr{\frac{Q}{1-\alpha}}^{\gamma^-}\pt \alpha}\\
&= \Dt\intOB{\frac{1}{\gamma^+-1} \lr{\frac{R}{\alpha}}^{\gamma^+}\alpha+\frac{1}{\gamma^--1}\lr{\frac{Q}{1-\alpha}}^{\gamma^-}(1-\alpha)}.
}
Multiplying the l.h.s. of the momentum equation \eqref{Div} by $-\Div\vu$ we therefore get
\eq{\label{ener_est}
&\sup_{t\in (0,T)}\intOB{\frac{1}{\gamma^+-1}\lr{\frac{R}{\alpha}}^{\gamma^+}\alpha+ \frac{1}{\gamma^--1}\lr{\frac{Q}{1-\alpha}}^{\gamma^-}(1-\alpha)}\\
&\qquad+ \intTO{ |\Div\vu|^2}\\
& \leq \intOB{\frac{1}{\gamma^+-1}\lr{\frac{R_0}{\alpha_0}}^{\gamma^+}\alpha_0+\frac{1}{\gamma^--1} \lr{\frac{Q_0}{1-\alpha_0}}^{\gamma^-}(1-\alpha_0)}.
}
On account of \eqref{equiv} and \eqref{R0Q0} the r.h.s. is bounded. Using \eqref{ener_est}, and \eqref{equiv} again, we obtain
\eqh{
\sup_{t\in (0,T)}\intO{Z^{\gamma^+}}= {\sup_{t\in (0,T)}}\intOB{Z^{\gamma^+}\alpha+ Z^{\gamma^+}(1-\alpha)}\\=
\sup_{t\in (0,T)}\intOB{Z^{\gamma^+}\alpha+ \lr{\frac{Q}{1-\alpha}}^{\gamma^-}(1-\alpha)}\leq C,
}
which provides the uniform estimate for $Z^{\gamma^+}$, as stated in \eqref{est_ener}.
The uniform estimate for $R^{\gamma^+}$ can be obtained using \eqref{RleqZ}, and then \eqref{equiv} provides the uniform estimate for $Q^{\gamma^-}$. 

In order to estimate the full gradient of $\vu$ we notice that
\eqh{\intTO{|\Grad\vu|^2}=\intTO{|\Div\vu|^2}+\intTO{|{\rm rot}\,\vu|^2},
}
the first term is bounded due to \eqref{ener_est}, and the second one is equal to $0$ on account of \eqref{rotu}.

In order to prove \eqref{est_bog} we multiply momentum equation \eqref{Div} by $P$ and integrating over time and space we get that
\eq{\label{P2}
\intTO{P^2}&=\intTO{\Div\vu\, P}+\intTO{\{P\} P}\\
&=\intTO{\Div\vu\, P}+\intT{\{P\}^2}.
}
The last term is bounded due to \eqref{est_ener}, and we use the Cauchy inequality to estimate
\eq{\label{divuP}
\left|\intTO{\Div\vu P}\right|&\leq \frac12\intTO{P^2}+\frac12\intTO{|\Div\vu|^2}\\
&\leq C+ \frac12\intTO{P^2}.
}
So, combining \eqref{divuP} with \eqref{P2} we obtain the uniform estimate for $Z^{2\gamma^+}$ as in \eqref{est_bog}. The rest of bounds from \eqref{est_bog} follows, as previously, from the relations between $Z$, $R$, and $Q$, see  \eqref{RleqZ}, and \eqref{equiv}. $\Box$

We can now use the above estimates in order to deduce that $R$ and $Q$ satisfy equations \eqref{SR} and \eqref{ST} in the renormalized sense.
\begin{lemma}\label{Lemma_ren}
Assume $X \in L^q((0,T)\times\T^d)$ with $q\geq 2$, and $\vu \in L^2(0,T;W^{1,2}(\T^d))$. Let $(X,\vu)$ solve 
$$\pt X+\Div(X\vu)=f~\text{in}~{\cal D}'((0,T)\times \T^d),$$
where $f\in L^p((0,T)\times\T^d)$ for some $p>1$, $p'\lr{\frac{q}{2}-1}\leq 1$.
Then $(X,\vu)$ is also a renormalized solution i.e. it solves 
\begin{equation}\label{renorent}
\partial_t b(X) + \Div(b(X)\vu) + \big(b'(X)X-b(X)\big) \Div\vu = f {b'(X)}~\text{in}~{\cal D}'((0,T)\times \T^d),
\end{equation}
where
\begin{equation}\label{regb}
b \in C^1([0,\infty)),\ | b'(s)|\leq Cs^\lambda,\ for\ s>1,\ where \ \lambda\leq \frac{q}{2}-1.
\end{equation}
\end{lemma}
The proof of this lemma is a consequence of the DiPerna-Lions theory \cite{DL} of renormalized solutions to the transport equation.

\bigskip

We now derive the equation satisfied by $Z$. If $R,Q,\vu$ are smooth, the evolution equation for $Z=Z(R,Q)$ can be deduced from the continuity equations for $R$ and for $Q$, and the formula \eqref{TZ}. However, we will use the the equivalency between the equations \eqref{SR}, \eqref{ST} and the evolution equation for $Z$ at the level of weak solutions, i.e. solutions with regularity specified in Lemma \ref{Lemma_est}. For such solutions we have the following result
\begin{lemma}\label{Lemma_equiv}
Let $\vu\in L^2(0,T; W^{1,2}(\T^d))$, $R\in L^{2\gamma^+}((0,T)\times \T^d)\cap L^\infty(0,T; L^{\gamma^+}(\T^d))$, 
$Q\in L^{2\gamma^-}((0,T)\times \T^d)\cap L^\infty(0,T; L^{\gamma^-}(\T^d))$,
and let $(R, Q,\vu)$ solve \eqref{SR} and \eqref{ST} in the sense of distributions, and {\eqref{Div} and \eqref{rotu} a.e. in $(0,T)\times\Omega$}. Then $Z$ defined by \eqref{TZ} belongs to $L^{2\gamma^+}((0,T)\times \T^d)\cap L^\infty(0,T; L^{\gamma^+}(\T^d))$ and it satisfies 
\eq{\label{SZ}
\pt Z+\Div(Z\vu)+\frac{(1-\gamma) (Z-R)Z}{ \gamma(Z-R)+R}\Div\vu=0,}
in the sense of distributions.
{Conversly}, let $(R,Z,\vu)$ solve \eqref{SR} and \eqref{SZ} in the sense of distributions. Then $Q$ defined by \eqref{TZ} satisfies \eqref{ST} in the sense of distributions. 
\end{lemma}
\pf The fact that $Z\in L^{2\gamma^+}((0,T)\times \T^d)\cap L^\infty(0,T; L^{\gamma^+}(\T^d))$ follows from Lemma \ref{Lemma_est}. Testing the equations \eqref{SR} and \eqref{ST} by $\xi_\eta(x-\cdot)$, where $\xi_\eta$ is a standard periodized mollifier, we obtain
\eq{\label{RQeta}
\pt R_\eta+\Div(R_\eta\vu)=r^1_\eta,\\
\pt Q_\eta+\Div(Q_\eta\vu)=r^2_\eta,
}
satisfied a.e. in $(0,T)\times\T^d$, where $a_\eta$ denotes $a\ast\xi_\eta$. From the Friedrichs commutator lemma we know that 
\eq{\label{r1r2}
r^1_\eta\to 0 \quad\mbox{in} \ L^{p_1}((0,T)\times \T^d),\\
r^2_\eta\to 0 \quad\mbox{in} \ L^{p_2}((0,T)\times \T^d)}
for $\frac{1}{p_1}=\frac{1}{2}+\frac{1}{2\gamma^+}$, $\frac{1}{p_2}=\frac{1}{2}+\frac{1}{2\gamma^-}$.
We now define $Z_\eta$ via
$$Q_\eta=\lr{1-\frac{R_\eta}{Z_\eta}}Z_\eta^\gamma,$$
and as previously we find that $R_\eta\leq Z_\eta$. Let us now apply $\partial_{Q_\eta}$, $\partial_{R_\eta}$ to both sides of the above formula,  we obtain respectively:
\[1=\gamma Z_\eta^{\gamma-1}\partial_{Q_\eta} Z_\eta-R_\eta(\gamma-1)Z_\eta^{\gamma-2}\partial_{Q_\eta} Z_\eta,\]
\[0=\gamma Z_\eta^{\gamma-1}\partial_{R_\eta} Z_\eta-Z_\eta^{\gamma-1}-R_\eta(\gamma-2)Z_\eta^{\gamma-1}\partial_{R_\eta} Z_\eta,
\]
therefore,
\eq{\label{pTR}
\partial_{Q_\eta} Z_\eta=\frac{1}{\gamma Z_\eta^{\gamma-1}-R_\eta(\gamma-1)Z_\eta^{\gamma-2}},\quad
\partial_{R_\eta} Z_\eta= \frac{Z_\eta^{\gamma-1}}{\gamma Z_\eta^{\gamma-1}-R_\eta(\gamma-1)Z_\eta^{\gamma-2}}.
}
Using the assumption $\gamma\leq1$  and inequality $R_\eta\leq Z_\eta$, we check that we can estimate the partial derivatives of $Z_\eta(R_\eta,Q_\eta)$ by the integrable function, more precisely
\eqh{
|\partial_{Q_\eta} Z_\eta|\leq\frac{Z_\eta^{1-\gamma}}{\gamma},\quad |\partial_{R_\eta} Z_\eta|\leq \frac{1}{\gamma}.
}
It means that $\partial_{Q_\eta} Z_\eta$ and $\partial_{R_\eta} Z_\eta$ are suitable test functions for  \eqref{RQeta}.
We now check that
\eqh{
\pt Z_\eta&=\partial_{Q_\eta}Z_\eta\pt Q_\eta+\partial_{R_\eta} Z_\eta\pt R_\eta\\
&= -\partial_{Q_\eta}Z_\eta\Div(Q_\eta\vu)-\partial_{R_\eta} Z_\eta\Div(R_\eta\vu)+r^1_\eta\partial_{R_\eta} Z_\eta+r^2_\eta\partial_{Q_\eta} Z_\eta\\
&= -[\partial_{Q_\eta}Z_\eta\Grad Q_\eta+\partial_{R_\eta} Z_\eta\Grad R_\eta]\cdot\vu-[Q_\eta\partial_{Q_\eta} Z_\eta +R_\eta\partial_{R_\eta} Z_\eta ]\Div\vu\\
&\quad+r^1_\eta\partial_{R_\eta} Z_\eta+r^2_\eta\partial_{Q_\eta} Z_\eta\\
&=-\Grad Z_\eta\cdot\vu-[Q_\eta\partial_{Q_\eta} Z_\eta +R_\eta\partial_{R_\eta} Z_\eta ]\Div\vu+r^1_\eta\partial_{R_\eta} Z_\eta+r^2_\eta\partial_{Q_\eta} Z_\eta,
}
therefore
\eq{\label{ZA}
\pt Z_\eta+\Div(Z_\eta\vu)+[Q_\eta\partial_{Q_\eta} Z_\eta +R_\eta\partial_{R_\eta} Z_\eta -Z_\eta]\Div\vu=r^1_\eta\partial_{R_\eta} Z_\eta+r^2_\eta\partial_{Q_\eta} Z_\eta.
}
Substituting \eqref{pTR} we compute
\eqh{
 &Q_\eta\partial_{Q_\eta} Z_\eta +R_\eta\partial_{R_\eta} Z_\eta-Z_\eta
 = \frac{Q_\eta+R_\eta Z_\eta^{\gamma -1}-Z_\eta(\gamma Z_\eta^{\gamma-1} - R_\eta(\gamma -1) Z_\eta^{\gamma-1})}{\gamma Z_\eta^{\gamma-1} - R_\eta(\gamma -1) Z_\eta^{\gamma-1}} \\
&= \frac{Z_\eta^\gamma - R_\eta Z_\eta^{\gamma-1} + R_\eta Z_\eta^{\gamma-1} - \gamma Z_\eta^\gamma + R_\eta(\gamma -1) Z_\eta^{\gamma -1}}{\gamma Z_\eta^{\gamma-1} - R_\eta(\gamma -1) Z_\eta^{\gamma-1}} =
 \frac{ (1-\gamma)(Z_\eta-R_\eta)Z_\eta}{\gamma(Z_\eta-R_\eta)+R_\eta}.
}
Using this, the equation for $Z_\eta$ can be written as
\eq{\label{Zs}
\pt Z_\eta+\Div(Z_\eta\vu)+\frac{(1-\gamma) (Z_\eta-R_\eta)Z_\eta}{ \gamma(Z_\eta-R_\eta)+R_\eta}\Div\vu=r^1_\eta\partial_{R_\eta} Z_\eta+r^2_\eta\partial_{Q_\eta} Z_\eta.
}
Note that since $\gamma <1$, we easily show that $\frac{(1-\gamma) (Z_\eta-R_\eta)Z_\eta}{ \gamma(Z_\eta-R_\eta)+R_\eta}\geq 0$, moreover $\frac{(1-\gamma) (Z_\eta-R_\eta)Z_\eta}{ \gamma(Z_\eta-R_\eta)+R_\eta}\leq \frac{1-\gamma}{ \gamma} (Z_\eta-R_\eta)\in L^{2\gamma^+}((0,T)\times\T^d)$. This allows us to let $\eta\to 0$ in the l.h.s. of \eqref{Zs}. The r.h.s. of \eqref{Zs} vanishes provided that $p_1\geq 1$ and $\frac{1}{p_2}+\frac{1-\gamma}{2\gamma^+}\leq 1$, which is fulfilled provided that $\gamma^+\geq1$.  

In order to recover  the equation for $Q$ from the equations for $R$ and $Z$, one derives the equation on $Q^{1/\gamma}$ first. It is easy to observe that \eqref{TZ} yields
$$
\left|\frac{\partial Q^{1/\gamma}}{\partial Z} \right|\leq 1 \quad\mbox{ and } \quad\left|\frac{\partial Q^{1/\gamma}}{\partial R} \right|\leq 1/\gamma.
$$
Therefore, the rigorous procedure involving mollifying  and passage to the limit $\eta\to 0$, described above, can be repeated for $Q$. In this manner we obtain a renormalized version of the equation for  $Q$ 
$$
\pt Q^{1/\gamma} + \Div (Q^{1/\gamma} \vu) + Q^{1/\gamma}\Div \vu=0.
$$
This finishes the proof. $\Box$

\bigskip

As a consequence of this Lemma and Lemma \ref{Lemma_ren} we have:
\begin{cor}
The couple $(Z,\vu)$ is a renormalized solution to \eqref{SZ}.
\end{cor}
\pf Indeed, for the assumptions of Lemma \ref{Lemma_ren} to be fulfilled we note that
$$\left|\frac{(1-\gamma) (Z-R)Z}{ \gamma(Z-R)+R}\Div\vu\right|\leq C Z|\Div\vu|,$$
and the r.h.s. is bounded in $L^1((0,T)\times\T^d)$ on account of \eqref{est_ener} and \eqref{est_bog}. $\Box$

\section{Sequential stability of solutions}\label{Sect:3}
The purpose of this section is to pass to the limit $n\to\infty$ in the sequence  $\{R_n, Q_n, Z_n, \vu_n\}_{n=1}^{\infty}$ and  to verify that the limit $(R,Q,Z,\vu)$ satisfies the system \eqref{S} in the weak sense. We  prove the following theorem
\begin{thm}\label{Thm_main2}
Let $T>0$. Assume that for any $n$ the quadruple 
 $(R_n, Q_n, Z_n,\vu_n)$ satisfies \eqref{TZ}--\eqref{S} with the initial conditions 
 \eqh{
 R_n \vert_{t=0} = R_{0,n}, \quad Q_n \vert_{t=0} = Q_{0,n}, \quad R_{0,n} \ge 0, \quad Q_{0,n} \ge 0,
 }
 with $Z_n \vert_{t=0} = Z_{0,n}$ satisfying \eqref{RelationInit},  and s.t.
 \eqh{
 R_{n,0}\to R_0\quad strongly \ in\ L^1(\T^d),\\
 Z_{n,0}\to Z_0\quad strongly \ in\ L^1(\T^d).
 }
Let the estimates from the Lemmas \ref{Lemma_posit} and \ref{Lemma_est} hold uniformly with respect to $n$. Then up to the subsequence
\eqh{
&R_n\to R\quad strongly\ in \ L^{2\gamma^+-\ep}((0,T)\times\T^d),\\
&Q_n\to Q\quad strongly \ in \  L^{2\gamma^--\ep}((0,T)\times\T^d),\\
&Z_n\to Z\quad strongly\ in \ L^{2\gamma^+-\ep}((0,T)\times\T^d),\\
&\vu_n\to \vu\quad weakly\ in\  L^2(0,T; H^1(\T^d)),
}
for any $\ep>0$. 
Moreover, $(R,Q,Z,\vu)$ satisfies \eqref{TZ}--\eqref{S} in the sense of distributions. 
\end{thm}
Passage to the limit $n\to\infty$ in the two first equations of system \eqref{S} requires at least weak convergence of the sequences $R_n, Q_n$ and $\vu_n$. This can be deduced directly from the a-priori estimates from Lemmas  \ref{Lemma_posit} and \ref{Lemma_est} using nowadays classical techniques (see, for example {\sc P.-L. Lions} \cite{Lions2} or {\sc E. Feireisl} \cite{EF2001}), and we skip this part.
The  core of the proof of Theorem \ref{Thm_main2} is to pass to the limit in the nonlinear term of the  momentum equation \eqref{Div}. Indeed, identification of the limit $\lim_{n\to\infty} p(R_n, Q_n)= p(R,Q)$ requires some sort of strong convergence of sequences $R_n$, and $Q_n$. Instead of proving the strong convergence of these sequences directly, we use  the equivalence between the system \eqref{S} and its reformulation in terms of $(R,Z,\vu)$, as stated in the Lemma \ref{Lemma_equiv}.
The proof of Theorem \eqref{Thm_main2} can be therefore reduced to the proof of compactness of the sequence $Z_n$ and justification that the limit quantities $R,\ Q,\ Z$ satisfy the relation \eqref{TZ}. We follow the strategy proposed by {\sc D. Bresch \& P.-E. Jabin} \cite{BJ}  (see also \cite{BJ_short}) in the context of compressible Navier-Stokes equations with the non-monotone pressure law.  As a byproduct of this approach, we obtain a compactness result for the sequence $R_n$, and using \eqref{TZ} the strong convergence of the sequence $Q_n$, concluding the proof of Theorem \ref{Thm_main2}. The rest of this section will be therefore devoted to the proof of the following result:
\begin{prop}\label{prop:rho}
Let $T>0$. Assume that $\{(R_n, Z_n,\vu_n)\}_{n=1}^{\infty}$ satisfies 
\begin{subequations}\label{S_Z}
\begin{align}
 &\pt R+ \Div(\vu R)=0,\label{SR_Z}\\ 
&\pt Z+\Div(Z\vu)+\frac{(1-\gamma) (Z-R)Z}{ \gamma(Z-R)+R}\Div\vu=0 ,\label{ST_Z}\\
 &{\Div\vu =Z^{\gamma^+}- \{Z^{\gamma^+}\}\label{Su_Z}},\\ 
&{\rm rot}\ \vu=0,\quad \intO{ \vu(t,x) }= 0,\label{Sur_Z}
\end{align}
\end{subequations}
 with the initial conditions \eqref{IcNew} satisfying \eqref{RelationInit}, and let the estimates from the Lemmas \ref{Lemma_posit} and \ref{Lemma_est} hold uniformly with respect to $n$. 
Then the sequences $\{R_n\}_{n=1}^\infty$ $\{Z_n\}_{n=1}^\infty$ are compact in $L^1((0,T)\times\T^d)$.
\end{prop}

\subsection{Preliminaries}\label{ssec:prelim}
In order to prove the strong convergence of $\{R_n,Z_n\}_{n=1}^{\infty}$ necessary to pass to the limit in the momentum equation, we will use the compactness criterion introduced in the context of Navier-Stokes equations in \cite{BJ}. 
First let us introduce the necessary notation.
We define the positive, bounded and symmetric function $\{K_h\}_{h>0}$ such that  
$$ K_h(x)=\frac{1}{(|x|+h)^{a}}$$ 
with $|x|=\sqrt {\sum_{i=1}^d x_i^2}$ for $|x|\leq 1/2$ with some $a>d$ and $K_h$ positive, independent of $h$ for $|x|\ge 2/3$, $K_h$ positive constant outside $B(0,3/4)$ so that $K_h\in{\cal C}^\infty(\T^d\backslash B(0,3/4))$ and it is a periodic function.
Further, we denote
\eq{\label{notationK}
\overline{K_h}(x) = \frac{K_h(x)}{\|K_h\|_{L^1(\T^d)}},
\qquad
\calK_{h_0}(x) = \int_{h_0}^1 \overline{K_h}(x) \frac{dh}{h}.
}
We  also use the following properties of the kernel $K_h$: 
\begin{equation}\label{ineqK}
K_h(x)= K_h(-x), \qquad |x||\nabla K_h(x)| \leq C K_h(x),
\end{equation}
for some constant $C>0$ independent of $h$ and
\eq{\label{Klogh}
\|{\cal K}_{h_0}\|_{L^1(\T^d)} \sim |\log h_0|.}

We recall the following compactness criterion, for the proof see \cite{Belgacem}, Lemma 3.1.
\begin{lemma}\label{lem:compact}
Let $\{X_n\}_{n=1}^\infty$ be a sequence of functions uniformly bounded in $L^p((0,T)\times\T^d)$
with $1\leq p<+\infty$. Assume that ${\cal K}_h$ is a sequence of positive, bounded functions s.t.
\begin{enumerate}
\item[i)] $\forall \eta>0$, \ $\sup_h \displaystyle \intO{{\cal K}_h(x)\vc{1}_{{\{x:\> |x|\geq\eta\}}}}<\infty$,
\item[ii)] $\|{\cal K}_h\|_{L^1(\T^d)}\to+\infty$\ as\ $h\to 0$.
\end{enumerate}
If $\{\pt X_n\}_{n=1}^\infty$ is uniformly bounded in $L^r([0,T],W^{-1,r}(\T^d))$ with 
$r\geq 1$ and  
$$
\underset{n}{\lim\sup} \left( \frac{1}{\|{\cal K}_h\|_{L^1}}  \int_0^T\!\!\!\iintO{
{\cal K}_h(x-y)|X_n(t,x)-X_n(t,y)|^p} \, \dt\right) \to 0, \quad \mbox{ as } h\to 0,
$$
then, $\{X_n\}_{n=1}^\infty$ is compact in $L^p([0,T]\times\T^d)$.
Conversely, if $\{X_n\}_{n=1}^\infty$ is compact in $L^p([0,T]\times\T^d)$, 
then the above $\lim\sup$ converges to $0$ as $h$ goes to $0$.
\end{lemma}
\subsection{Propagation of oscillations }
Having  transport equations for $R$ and $Z$ together with necessary a-priori bounds, our next goal is to derive the equations for perturbations of both of these quantities. 
Perturbations are described by the evolution of $|R(t,x)-R(t,y)|$ and $|Z(t,x)-Z(t,y)|$, respectively, for any couple of points $x,y\in\T^d$. To obtain them we first subtract the equations for $R(t,x)$ and $R(t,y)$  
\eqh{
&\pt (R_x-R_y) + \Div_x (\vu_x\lr{R_x-R_y}) +
\Div_y (\vu_y\lr{R_x-R_y})  \\
&=\frac 12 (\Div_x \vu_x + \Div_y \vu_y) \lr{R_x-R_y}  - \frac 12 (\Div_x \vu_x-\Div_y \vu_y)(R_x+R_y),
}
where we denoted $R_x=R(t,x), \ R_y=R(t,y)$.
Multiplying this equation by the sign of their difference $s_R=sign(R_x-R_y)$ we get
\eq{\label{eqR}
&\pt |R_x-R_y| + \Div_x (\vu_x|R_x-R_y|) +
\Div_y (\vu_y|R_x-R_y|)  \\
&=\frac 12 (\Div_x \vu_x + \Div_y \vu_y) |R_x-R_y|  - \frac 12 (\Div_x \vu_x-\Div_y \vu_y)(R_x+R_y)s_R.
}
By the similar token, we obtain the equation for $|Z_x-Z_y|$, namely
\eq{\label{eqZ}
&\pt|Z_x-Z_y|+
\Div_{x}(\vu_x|Z_x-Z_y|)+ \Div_{y}(\vu_y|Z_x-Z_y|)\\
&=\frac 12 (\Div_x \vu_x + \Div_y \vu_y) |Z_x-Z_y| - \frac 12 (\Div_x \vu_x-\Div_y \vu_y)(Z_x+Z_y)s_Z\\
&\quad-\left[
\frac{(1-\gamma) (Z_x-R_x)Z_x}{ \gamma(Z_x-R_x)+R_x} \Div_x \vu_x-
\frac{(1-\gamma) (Z_y-R_y)Z_y}{ \gamma(Z_y-R_y)+R_y}\Div_y\vu_y
\right] s_Z.
}
We  now multiply \eqref{eqR} and \eqref{eqZ} by $w_x+w_y$, where $w_x=w(t,x)$ denotes the solution to the transport equation
\eq{\label{def_w}
\left\{\begin{array}{l}
\pt w+\vu\cdot\Grad w+ \theta {\cal D} w=0,\\
w(0,x)=1,
\end{array}\right.
}
where $\theta$ is a constant parameter that will be chosen later on and ${\cal D}$ depending on $\vu$ and $Z$.
Our next step is to write equation for 
\eqh{
S(t):=\iintO{ K_h(x-y)O_{x-y} (w_x+w_y)},
}
where $O_{x-y}=|R_x-R_y|+ |Z_x-Z_y|$, we have
\eq{\label{Rt0}
\Dt S(t)&= \iintO{ \Grad K_h(x-y)(\vu_x-\vu_y)O_{x-y} (w_x+w_y)}\\
&\quad+\frac 12 \iintO{K_h(x-y)(\Div_x \vu_x + \Div_y \vu_y)O_{x-y} (w_x+w_y)}\\
&\quad - \frac 12 \iintO{K_h(x-y)(\Div_x \vu_x-\Div_y \vu_y)[(R_x+R_y)s_R+(Z_x+Z_y)s_Z](w_x+w_y)}\\
&\quad-\iintO{K_h(x-y)\left[
\frac{(1-\gamma) (Z_x-R_x)Z_x}{ \gamma(Z_x-R_x)+R_x} \Div_x \vu_x-
\frac{(1-\gamma) (Z_y-R_y)Z_y}{ \gamma(Z_y-R_y)+R_y}\Div_y\vu_y
\right] s_Z(w_x+w_y)}\\
&\quad+ \iintO{K_h(x-y)O_{x-y} \lr{\pt w_x+\vu_x\cdot\Grad w_x}}\\
&\quad+ \iintO{K_h(x-y)O_{x-y} \lr{\pt w_y+\vu_y\cdot\Grad w_y}}
}
Using the symmetry of $K_h(x-y)$, $O_{x-y}$, and the symmetry of the second, third and fourth integrals on the r.h.s. we obtain
 \eq{\label{Rt1}
&\Dt S(t)= \iintO{ \Grad K_h(x-y)(\vu_x-\vu_y)O_{x-y} (w_x+w_y)}\\
&\quad+ \iintO{K_h(x-y)(\Div_x \vu_x + \Div_y \vu_y)O_{x-y} w_x}\\
&\quad - \iintO{K_h(x-y)(\Div_x \vu_x-\Div_y \vu_y)[(R_x+R_y)s_R+(Z_x+Z_y)s_Z]w_x}\\
&\quad-2\iintO{K_h(x-y)\left[
\frac{(1-\gamma) (Z_x-R_x)Z_x}{ \gamma(Z_x-R_x)+R_x} \Div_x \vu_x-
\frac{(1-\gamma) (Z_y-R_y)Z_y}{ \gamma(Z_y-R_y)+R_y}\Div_y\vu_y
\right] s_Zw_x}\\
&\quad+2 \iintO{K_h(x-y)O_{x-y} \lr{\pt w_x+\vu_x\cdot\Grad w_x}}.
}
Finally, writing $(\Div_x \vu_x + \Div_y \vu_y)=-(\Div_x \vu_x - \Div_y \vu_y)+2\Div_x \vu_x$ and combining the second and the third term on the r.h.s. we obtain
 \eq{\label{Rt2}
&\Dt S(t)= \iintO{ \Grad K_h(x-y)(\vu_x-\vu_y)O_{x-y} (w_x+w_y)}\\
&\quad - \iintO{K_h(x-y)(\Div_x \vu_x-\Div_y \vu_y)[R_xs_R+Z_x s_Z]w_x}\\
&\quad-2\iintO{K_h(x-y)\left[
\frac{(1-\gamma) (Z_x-R_x)Z_x}{ \gamma(Z_x-R_x)+R_x} \Div_x \vu_x-
\frac{(1-\gamma) (Z_y-R_y)Z_y}{ \gamma(Z_y-R_y)+R_y}\Div_y\vu_y
\right] s_Zw_x}\\
&\quad+2 \iintO{K_h(x-y)O_{x-y} \lr{\pt w_x+\vu_x\cdot\Grad w_x+\Div_x\vu_x w_x}}\\
&=I_1+I_2+I_3+I_4.
}
We now estimate each term in \eqref{Rt2}, separately.

\medskip

\noindent{\bf Estimate of $I_1$.} Recall that $K_h$ satisfies \eqref{ineqK}, 
we also know that
\eq{\label{est_Du}
|\vu_x-\vu_y|\leq C|x-y|(D_{|x-y|}\vu_x+D_{|x-y|}\vu_y),\quad\mbox{where}\quad D_h \vu_x=\frac{1}{h}\int_{|z|\leq h}\frac{|\Grad \vu_{x+z}|}{|z|^{d-1}}\,{\rm d}z.}
Recall that $D_h \vu  \leq M |\nabla \vu|$, where $M$  denotes the maximal operator:
$$
M f(x) = \sup_{r\leq 1} \frac{1}{|B(0,r)|} \int_{B(0,r)} f(x+z)\,dz.
$$
For the proof of this fact we refer the reader to \cite[Lemma 3.1 and Eq. 3.3]{Jabin2010}.
Combining estimate \eqref{est_Du} with \eqref{ineqK}, we have
\eqh{
I_1&=\iintO{ \Grad K_h(x-y)(\vu_x-\vu_y)O_{x-y} (w_x+w_y)}\\
&\leq C\iintO{K_h(x-y)(D_{|x-y|}\vu_x+D_{|x-y|}\vu_y)O_{x-y} w_x}.
}
Next, writing
$$
D_{|x-y|}\vu_x+D_{|x-y|}\vu_y =  D_{|x-y|}\vu_y-D_{|x-y|}\vu_x + 2D_{|x-y|}\vu_x,
$$ 
and estimating $D_{|x-y|}\vu(x)$ by the Maximal operator $M|\Grad\vu|(x)$ we get
\eq{\label{est_I1}
I_1\leq &C\iintO{K_h(x-y)(D_{|x-y|}\vu_y-D_{|x-y|}\vu_x)O_{x-y} w_x}\\
&+C\iintO{K_h(x-y)M|\Grad\vu_x|O_{x-y} w_x}.
}

\medskip

\noindent{\bf Estimate of $I_2$.}  {Due to \eqref{Su_Z}, we can write that}
\eqh{
 -\lr{ \Div_x\vu_x-\Div_y\vu_y}[R_xs_R+Z_x s_Z]
 =- \lr{ Z_x^{\gamma^+}-Z_y^{\gamma^+}}R_xs_R- \lr{ Z_x^{\gamma^+}-Z_y^{\gamma^+}}Z_x s_Z.
 }
 Note that the last term is always nonpositive, hence the contribution to $I_1$ coming from this term can be moved to the l.h.s. of \eqref{Rt2}. Concerning the first term, it can only be estimated from above. Using $R_x=\alpha_x Z_x$ we get
 \eq{\label{est_I2}
I_2= &- \iintO{K_h(x-y)(\Div_x \vu_x-\Div_y \vu_y)[R_xs_R+Z_x s_Z]w_x}\\
\leq&- \iintO{K_h(x-y)|Z_x^{\gamma^+}-Z_y^{\gamma^+}|Z_xw_x}\\
&+
\iintO{K_h(x-y)|Z_x^{\gamma^+}-Z_y^{\gamma^+}|\alpha_xZ_xw_x}\\
=&\iintO{K_h(x-y)|Z_x^{\gamma^+}-Z_y^{\gamma^+}|(\alpha_x-1)Z_xw_x}.
 }
 Since $\alpha_x\leq1$, the integral $I_2$ is nonpositive and can be moved to the left hand side of \eqref{Rt2}.
\medskip

\noindent{\bf {Estimate of $I_3$}.}
 The estimate of this term is most lengthy and requires splitting the term in a big bracket several times. We first note that
\eq{\label{re1}
&-\left[
\frac{(1-\gamma) (Z_x-R_x)Z_x}{ \gamma(Z_x-R_x)+R_x} \Div_x \vu_x-
\frac{(1-\gamma) (Z_y-R_y)Z_y}{ \gamma(Z_y-R_y)+R_y}\Div_y\vu_y
\right] s_Z\\
&=(\gamma-1)\left[
\frac{ (Z_x-R_x)Z_x^{\gamma^++1}}{ \gamma(Z_x-R_x)+R_x}-
\frac{(Z_y-R_y)Z_y^{\gamma^++1}}{ \gamma(Z_y-R_y)+R_y}
\right] s_Z\\
&-(\gamma-1)\left[ \frac{ (Z_x-R_x)Z_x}{ \gamma(Z_x-R_x)+R_x}-
\frac{(Z_y-R_y)Z_y}{ \gamma(Z_y-R_y)+R_y} 
\right] s_Z \{Z^{\gamma^+}\}\\
&=(\gamma-1)\frac{(Z_y-R_y)}{ \gamma(Z_y-R_y)+R_y}\lr{Z_x^{\gamma^++1}-Z_y^{\gamma^++1}} s_Z\\
&+(\gamma-1)
\left[
\frac{ (Z_x-R_x)}{ \gamma(Z_x-R_x)+R_x}-
\frac{(Z_y-R_y)}{ \gamma(Z_y-R_y)+R_y}
\right] Z_x^{\gamma^++1} s_Z\\
&-(\gamma-1)\frac{(Z_y-R_y)}{ \gamma(Z_y-R_y)+R_y}\lr{Z_x-Z_y} s_Z \{Z^{\gamma^+}\} \\
&-(\gamma-1)
\left[
\frac{ (Z_x-R_x)}{ \gamma(Z_x-R_x)+R_x}-
\frac{(Z_y-R_y)}{ \gamma(Z_y-R_y)+R_y}
\right] Z_x\{Z^{\gamma^+} \}s_Z,
}
{where we used \eqref{Su_Z} to substitute} for $\Div\vu$, and the observation that $\{Z_x^{\gamma^+}\}=\{Z_y^{\gamma^+}\}$.

Observe that
\eq{\label{sign_frac}
\frac{(Z_y-R_y)}{ \gamma(Z_y-R_y)+R_y}=\frac{(1-\alpha_y)}{ \gamma(1-\alpha_y)+\alpha_y}\ \in  \lr{0,1/{\gamma}},
}
 therefore, the first term on the r.h.s.  of \eqref{re1} is non-positive for $\gamma<1$ and can be eventually considered on the l.h.s. of \eqref{Rt2}. For the second term on the r.h.s.  of \eqref{re1} we continue to write
 \eqh{
 &(\gamma-1)
\left[
\frac{ (Z_x-R_x)}{ \gamma(Z_x-R_x)+R_x}-
\frac{(Z_y-R_y)}{ \gamma(Z_y-R_y)+R_y}
\right] Z_x^{\gamma^++1} s_Z\\
&=(\gamma-1)\frac{R_y Z_x-R_xZ_y}{\lr{ \gamma(Z_x-R_x)+R_x}\lr{ \gamma(Z_y-R_y)+R_y}}Z_x^{\gamma^++1} s_Z\\
&=(\gamma-1)\frac{R_y (Z_x-Z_y)}{\lr{ \gamma(Z_x-R_x)+R_x}\lr{ \gamma(Z_y-R_y)+R_y}}Z_x^{\gamma^++1} s_Z\\
&+(\gamma-1)\frac{(R_y -R_x)Z_y}{\lr{ \gamma(Z_x-R_x)+R_x}\lr{ \gamma(Z_y-R_y)+R_y}}Z_x^{\gamma^++1} s_Z,
}
 again, the first term has a good sign, while the second one can be transformed to
 \eqh{
 &(\gamma-1)\frac{(R_y -R_x)Z_y}{\lr{ \gamma(Z_x-R_x)+R_x}\lr{ \gamma(Z_y-R_y)+R_y}}Z_x^{\gamma^++1} s_Z\\
 &=(\gamma-1)\frac{(R_y -R_x)}{\lr{ \gamma(1-\alpha_x)+\alpha_x}\lr{ \gamma(1-\alpha_y)+\alpha_y}}Z_x^{\gamma^+} s_Z\\
& \leq\frac{1-\gamma}{\gamma^2}|R_x-R_y| Z_x^{\gamma^+}.
 }
The third and fourth terms on the r.h.s. of \eqref{re1} can be treated similarly and estimated by 
$\frac{1-\gamma}{\gamma}|Z_x-Z_y|\{Z^{\gamma^+} \}$, and  by $\frac{1-\gamma}{\gamma^2}O_{x-y}\{Z^{\gamma^+} \}$, respectively.

 Putting all the terms together, we conclude that
 \eq{\label{est_I3}
&I_3=-2\iintO{K_h(x-y)\left[
\frac{(1-\gamma) (Z_x-R_x)Z_x}{ \gamma(Z_x-R_x)+R_x} \Div_x \vu_x-
\frac{(1-\gamma) (Z_y-R_y)Z_y}{ \gamma(Z_y-R_y)+R_y}\Div_y\vu_y
\right] s_Zw_x}\\
&\leq -2(1-\gamma)\iintO{K_h(x-y)\frac{(Z_y-R_y)}{ \gamma(Z_y-R_y)+R_y}\left|Z_x^{\gamma^++1}-Z_y^{\gamma^++1}\right|w_x}\\
&\quad-2(1-\gamma)\iintO{K_h(x-y)\frac{R_y |Z_x-Z_y|}{\lr{ \gamma(Z_x-R_x)+R_x}\lr{ \gamma(Z_y-R_y)+R_y}}Z_x^{\gamma^++1}w_x}\\
&\quad+C\iintO{K_h(x-y)|R_x-R_y| Z_x^{\gamma^+}w_x}+C\{Z^{\gamma^+} \}\iintO{K_h(x-y)O_{x-y}w_x},
 }
 with some constant $C$ depending on $\gamma$.
\medskip

\noindent{\bf Estimate of $I_4$.}
  For the last term in \eqref{Rt2} we simply use the definition of the measure $w$ from \eqref{def_w}, we therefore get
  \eq{\label{est_I4}
  I_4=&2 \iintO{K_h(x-y)O_{x-y} \lr{\pt w_x+\vu_x\cdot\Grad w_x+\Div_x\vu_x w_x}}\\
  =&2 \iintO{K_h(x-y)O_{x-y} (\Div_x\vu_x-\theta {\cal D}_x) w_x}.
  }
  As a conclusion, we obtain from \eqref{Rt2}, using estimates \eqref{est_I1}, \eqref{est_I2}, \eqref{est_I3}, and \eqref{est_I4} that
  \eq{\label{mainS}
  \Dt S(t)&+\iintO{K_h(x-y)|Z_x^{\gamma^+}-Z_y^{\gamma^+}|(1-\alpha_x)Z_xw_x}\\
  &+2(1-\gamma)\iintO{K_h(x-y)\frac{(Z_y-R_y)}{ \gamma(Z_y-R_y)+R_y}\left|Z_x^{\gamma^++1}-Z_y^{\gamma^++1}\right|w_x}\\
  &+2(1-\gamma)\iintO{K_h(x-y)\frac{R_y |Z_x-Z_y|}{\lr{ \gamma(Z_x-R_x)+R_x}\lr{ \gamma(Z_y-R_y)+R_y}}Z_x^{\gamma^++1}w_x}\\
  \leq
  &C\iintO{K_h(x-y)(D_{|x-y|}\vu_y-D_{|x-y|}\vu_x)O_{x-y} w_x}\\
&+C\iintO{K_h(x-y)M|\Grad\vu_x|O_{x-y} w_x}\\
&+C\iintO{K_h(x-y)|R_x-R_y| Z_x^{\gamma^+}w_x}+C\{Z^{\gamma^+} \}\iintO{K_h(x-y)O_{x-y}w_x}\\
&+2 \iintO{K_h(x-y)O_{x-y} (\Div_x\vu_x -\theta {\cal D}_x)w_x}.
  }
Note that the only assumptions are: $\gamma<1$, i.e. $\gamma^+<\gamma^-$, and $\alpha\in[0,1]$.

\subsection{Compactness criterion with weights}
At this point of the proof it is convenient to chose the ${\cal D}$ from the definition of the weight $w_x$ \eqref{def_w}. Taking for example
\eq{\label{def_D}
{\cal D}=M|\Grad\vu|+|\Div\vu|+Z^{\gamma^+}+ \{Z^{\gamma^+}\},}
and $\theta$ sufficiently large we obtain from \eqref{mainS} 
\eqh{\Dt S(t)\leq C\iintO{K_h(x-y)(D_{|x-y|}\vu_y-D_{|x-y|}\vu_x)O_{x-y} w_x}.}
Note that, thanks to uniform estimates from Lemmas \ref{Lemma_posit} and \ref{Lemma_est},
we have that ${\cal D}$ defined in \eqref{def_D} is uniformly bounded in $L^2((0,T)\times\T^d)$ (note that $\|M|\Grad\vu|\|_{L^2((0,T)\times\T^d)}\leq C\|\Grad\vu\|_{L^2((0,T)\times\T^d)}$). This allows us to deduce the following properties of the weight $w$. 

\begin{prop}[Proposition 7.2 in \cite{BJ}]
\label{prop:w}
Assume that ${\cal D}$ defined in \eqref{def_D} is given and that it is  bounded in $L^2(0,T\times\T^d)$.
Then, there exists a weight $w$ solving  \eqref{def_w}. Moreover, we have
\begin{itemize}
\item[(i)] For any $(t,x)\in(0,T)\times T^d$, $0\leq w(t,x) \leq 1$.
\item[(ii)] If we assume moreover that the pair $(X,\vu)$ is a solution to the continuity equation:
\eqh{\pt X+\Div(X\vu)=0,}
 and $X$ is bounded in $L^2((0,T)\times\T^d)$, there exists $C\geq 0$, such that
\begin{equation}\label{boundlogw}
\intO{ X |\log w| } \leq C \theta.
\end{equation}
\end{itemize}
\end{prop}
Let us now introduce
\eq{\label{Sfini}
S_{h_0}(t)=\int_{h_0}^1 \frac{S(t)}{\|K_h\|_{L^1}}\,\frac{dh}{h}=\int_{h_0}^1\int_{\T^{2d}} \overline{K_h}(x-y)O_{x-y} (w_x+w_y)\, \dx\, \dy\,\frac{{\rm d}h}{h}.}
Recalling the notation \eqref{notationK}, and changing the variables $z=x-z$ we get from \eqref{Sfini} that
\eqh{\Dt S_{h_0}(t)\leq &C \int_{h_0}^1 \iintO{\overline{K_h}(x-y) \left|D_{|x-y|}\vu(y)-D_{|x-y|}\vu(x)\right|O_{x-y} w_x}\,\frac{{\rm d}h}{h}\\
\leq &C \int_{h_0}^1 \int_{\T^{d}}\overline{K_h}(z) \|D_{|z|}\vu(\cdot)-D_{|z|}\vu(\cdot+z)\|_{L^2} \,{\rm d}z\,\frac{{\rm d}h}{h}
,}
where to get the last inequality we have used the Cauchy-Schwarz inequality, $L^2$ integrability of $R$ and $Z$ together with boundedness of the weight $w$. Integrating with respect to time we obtain
\eq{\label{Rt0}
S_{h_0}(t)-S_{h_0}(0)\leq  &C\int_0^t\int_{h_0}^1 \int_{\T^{d}}\overline{K_h}(z) \|D_{|z|}\vu(\cdot)-D_{|z|}\vu(\cdot+z)\|_{L^2} \,{\rm d}z\,\frac{{\rm d}h}{h}\,{\rm d}s.
}
We now use the following lemma from \cite{BJ}.
\begin{lemma}[Lemma 6.3 in \cite{BJ}]
{There exists $C>0$ such that} for any $\vu\in H^1(\T^d)$,
\begin{equation}\label{ineq:lem63}
\int_{h_0}^1 \int_{\T^{d}} \overline{K_h}(z) \|D_{|z|}\vu - D_{|z|}\vu(\cdot+z)\|_{L^2(\T^d)} \,{\rm d}z \frac{{\rm d}h}{h}
\leq C |\log h_0|^{1/2} \|\vu\|_{H^1(\T^d)}.
\end{equation}
\end{lemma}
With this at hand \eqref{Rt0} gives
\eqh{
S_{h_0}(t)-S_{h_0}(0)\leq  C|\log(h_0)|^{1/2} \int_0^t \|\vu(s)\|_{H^1(\T^d)}\, {\rm d}s.
}
From this, boundedness of $\vu$ in $L^2(0,T; H^1(\T^d))$ we get that
\eq{\label{Sg}
S_{h_0}(t)\leq C|\log(h_0)|^{1/2}+ S_{h_0}(0).}
Using the reversed Lemma \ref{lem:compact} and strong convergence of initial data we get that  
\eqh{
&\limsup_n\frac{S_{h_0}(0)}{|\log h_0|}\\
&=\limsup_n\lr{\frac{\iintO{ \calK_{h_0}(x-y)\Big[ \left|Z_{0,n}(x)-Z_{0,n}(y)\right|+ |R_{0,n}(x)-R_{0,n}(y)|\Big] }}{|\log h_0|}}\to 0,
}
as $h_0\to 0$. Therefore, \eqref{Sg} yields that
\eq{\label{Sw}
\sup_{t\in[0,T]}\frac{S_{h_0}}{|\log h_0|}\to 0,\quad \mbox{as}\quad h_0\to0.
}

\subsection{Removal of the weights}
We now want to remove the weights from \eqref{Sw} to prove that
\eq{\label{without_w}
\limsup_n\lr{\frac{1}{|\log h_0|}\iintO{ \calK_{h_0}(x-y)\Big[ \left|(Z_x)_n-(Z_y)_n\right|+ |(R_x)_n-(R_y)_n|\Big] }}\to 0
}
as $h_0\to 0$, 
while for the moment we only have only convergence with the weights, see \eqref{Sfini} and \eqref{Sw}.
{We present a formal
argument leading to the estimate \eqref{boundlogw} for $X=Z_n$}. It is clear that estimate \eqref{boundlogw} holds for $X=R_n$, we will explain why one can expect the same estimate for $X=Z_n$.
Recall that from \eqref{def_w} know that the weights satisfy the equation
\eq{\label{logw}
\pt|\log w_n|+\vu_n\cdot \Grad |\log w_n|=\theta {\cal D}_n,}
with
\eqh{
{\cal D}_n=M|\Grad\vu_n|+|\Div\vu_n|+Z_n^{\gamma^+}+ \{Z_n^{\gamma^+}\},}
Multiplying \eqref{logw} by $R_n$ and using the fact that $R_n$ satisfies the continuity equation we get
\eq{\label{cwR}
\Dt\intO{R_n|\log w_n|}= \theta\intO{R_n{\cal D}_n}.}
Mimicking this procedure for $Z_n$ satisfying equation \eqref{ST_Z} we have
\eq{\label{Zpre}
\Dt\intO{Z_n|\log w_n|}=\theta\intO{Z_n{\cal D}_n}- \intO{\frac{(1-\gamma) (Z_n-R_n)Z_n}{ \gamma(Z_n-R_n)+R_n}\Div\vu_n|\log w_n|}.
}
The integrability of the terms $R_n{\cal D}_n$ and $Z_n{\cal D}_n$ on r.h.s. of {\eqref{cwR} and \eqref{Zpre}, respectively}, follows from the fact that $R_n{\cal D}_n\leq Z_n{\cal D}_n$ is bounded in $L^1((0,T)\times\T^d)$. Indeed, on account of Lemma \ref{Lemma_est} and the assumption $\gamma^+\geq1$, we have $Z_n^{\gamma^++1} \in L^{1}((0,T)\times \T^d)$. Therefore from \eqref{cwR} it follows that
\eq{\label{lwR}
\intO{R_n|\log w_n|}\leq C\theta.}
To deduce the same inequality for $Z_n$, we rewrite the term $\Div\vu_n$  in \eqref{Zpre} using the formula \eqref{Su_Z} 
\eqh{
\Dt\intO{Z_n|\log w_n|}+\intO{\frac{(1-\gamma) (Z_n-R_n)Z_n}{ \gamma(Z_n-R_n)+R_n}Z_n^{\gamma^+}|\log w_n|}\\
=\theta\intO{Z_n{\cal D}_n}- \intO{\frac{(1-\gamma) (Z_n-R_n)Z_n}{ \gamma(Z_n-R_n)+R_n}\{Z_n^{\gamma^+}\}|\log w_n|}
.} 
 Recalling \eqref{sign_frac}, the second term on the l.h.s. is nonnegative, while the second term on the r.h.s. can be bounded, therefore we have
 \eq{\label{cwZ}
\Dt\intO{Z_n|\log w_n|}
\leq \theta\intO{Z_n{\cal D}_n}+\frac{1-\gamma}{\gamma}\{Z_n^{\gamma^+}\}\intO{Z_n|\log w_n|}
.} 
Applying the Gronwall Lemma we can show that
\eqh{
&\intO{Z_n|\log w_n|(t)}\\
&\leq \exp\left\{\frac{1-\gamma}{\gamma}\int_0^t\{Z_n^{\gamma^+}(s)\}\,{\rm d} s\right\}
\lr{\intO{Z_n|\log w_n|(0)}+ \theta\int_0^t\intO{Z_n{\cal D}_n(t)}\,{\rm d} s}.
}
So, due to \eqref{def_w}, Lemma \ref{Lemma_est} and the discussion above we get that
\eq{\label{lwZ}
\intO{Z_n|\log w_n|}\leq C\theta.}
The above reasoning may be made rigorous by following the proof of the Proposition \ref{prop:w}  presented in \cite[Proposition 7.2]{BJ} with minor changes due to presence of additional terms in the equation for $Z_n$.

Let $\eta<1$. We define $\omega_\eta = \{x: w \leq \eta\}$ and denote by 
$\omega_\eta^c$ its complementary.
We have 
\begin{align}
\iintO{\calK_{h_0}(x-y)O_{x-y}} 
&= \int_{h_0}^1 \iintO{ \overline{K_h}(x-y)O_{x-y}}\frac{{\rm d}h}{h} = A_1 + A_2, \label{eqI12}
\end{align}
with
\begin{align*}
A_1 & = \int_{h_0}^1 \int_{\{x\in\omega_\eta^c\}\cup\{y\in\omega_\eta^c\}} \overline{K_h}(x-y) O_{x-y}\dx\,\dy\frac{{\rm d}h}{h} \leq \frac{1}{\eta} S_{h_0},
\end{align*}
due to \eqref{Sfini}, 
and
\begin{align*}
A_2 & =\int_{h_0}^1 \int_{\{x\in\omega_\eta\}\cap\{y\in\omega_\eta\}} \overline{K_h}(x-y) O_{x-y}\,\dx\,\dy\frac{{\rm d}h}{h}  \\
& \leq 2 \int_{h_0}^1 \iintO{\overline{K_h}(x-y) [R(x)+Z(x)]\frac{|\log w(x)|}{|\log \eta|}}\frac{{\rm d}h}{h},
\end{align*}
where we use the symmetry of $K_h$ and the fact that, by definition, 
for $\eta<1$, $|\log w(x)|\geq |\log \eta|$ for all $x\in\omega_\eta$.
Changing the variables $z=y-x$, and recalling $\|\overline{K_h}\|_{L^1}=1$,
we get
$$
A_2 \leq \frac{2 |\log h_0|}{|\log \eta|} \intO{[ R(x)+Z(x)]|\log w(x)| },
$$
where the integral on the r.h.s. is bounded due to \eqref{lwR} and \eqref{lwZ}.

Summarizing the estimates of $A_1$, $A_2$, we obtain
$$
\intO{ \calK_{h_0}(x-y)\Big[ \left|(Z_x)_n-(Z_y)_n\right|+ |(R_x)_n-(R_y)_n|\Big] } 
\leq \frac{1}{\eta} S_{h_0} + C\frac{|\log h_0|}{|\log \eta|}.
$$
 Applying \eqref{Sg}, we obtain
$$
\iintO{ \calK_{h_0}(x-y)\Big[ \left|(Z_x)_n-(Z_y)_n\right|+ |(R_x)_n-(R_y)_n|\Big] }
\leq C\frac{|\log h_0|^{1/2}+1}{\eta}  + C\frac{|\log h_0|}{|\log \eta|},
$$
Next, due to \eqref{Klogh}, we have
\eqh{
\frac{1}{\|\calK_{h_0}\|_{L^1}}\iintO{ \calK_{h_0}(x-y)\Big[ \left|(Z_x)_n-(Z_y)_n\right|+ |(R_x)_n-(R_y)_n|\Big]}
\\
\leq C\frac{|\log h_0|^{1/2}+1}{\eta|\log h_0|}  + \frac{C}{|\log \eta|},
}
so, choosing for example $\eta=|\log h_0|^{-1/4}$, and letting $h_0\to 0$, we show \eqref{without_w}, which proves the compactness criterion from Lemma \ref{lem:compact}.
\subsection{Concluding remarks}
From the previous section it follows that the sequences $R_n$ and $Z_n$ converge strongly in $L^1((0,T)\times\T^d)$ to $R,\ Z$, respectively. Using the uniform bounds from Lemma \ref{Lemma_est}, we therefore deduce the strong convergence of both of these sequences in $L^{2\gamma^+-\ep}((0<T)\times\T^d)$, for any $\ep>0$. Using the equivalence relation \eqref{TZ}, and the uniform estimates on $Q_n$, we deduce that $Q_n\to Q\quad strongly \ in \  L^{2\gamma^--\ep}((0<T)\times\T^d)$. With this information at hand, it is possible to pass to the limit in all equations of system  \eqref{TZ}--\eqref{S}. The proof of Theorem \ref{Thm_main2} is therefore complete. $\Box$

\section{Existence of approximate solutions}\label{Sect:4}
Here we prove the existence of solutions to certain approximation of system \eqref{S} involving parameter $k$.  Let us denote  ${\cal T}_k(\cdot)$ the truncation operator $ {\cal T}_k: \R_+ \to \R_+$ for $k\in \R_+$ such that
\eqh{
 {\cal T}_k(t)= t \mbox{ for } t \le k\quad  \mbox{ and }\quad  {\cal T}_k(t)=k \mbox{ for }
    t\geq k.
}
We consider the following approximate system for $(t,x)\in (0,T)\times \T^d$:
\begin{subequations}\label{S-appr}
\begin{align}
  &\d_t R+ \Div(\vu R) 
  =0,\label{S-appr1}\\
  &\d_t Q + \Div(\vu Q) 
  =0 ,\label{S-appr2} \\
  &Q = \lr{1-\frac{R}{Z}} Z^\gamma,\quad R\leq Z,\label{S-appr3}\\
   &{ \div \vu = ({\cal T}_k(Z))^{\gamma_+} -\{({\cal T}_k(Z))^{\gamma_+}\}} \label{S-appr4},\\
  &{\rm rot} \, \vu = \vc{0},\label{S-appr5}\\
  &\intO{\vu(t,x)}={0}.\label{S-appr6}
\end{align}
\end{subequations}
Note that combination of \eqref{S-appr5} with \eqref{S-appr6} defines $\vu$ as a potential flow, i.e. there exists $\phi(t,x)$ such that $\vu=\nabla \phi$, equivalently $\Delta \phi =\div \vu$. 
The main result of this section is existence of solutions to the approximate system \eqref{S-appr} supplemented by the initial conditions
 \eq{\label{ini_ep}
 R_k \vert_{t=0} = R_{0,k}, \quad Q_k \vert_{t=0} = Q_{0,k},
 }
 and $ Z_k \vert_{t=0} = Z_{0,k}$ defined by \eqref{RelationInit}.

\begin{thm}\label{approx}
Let the initial conditions \eqref{ini_ep} be  such that
\eqh{0< R_{0,k},\, Q_{0,k}<\infty \ \text{a.e. in } \T^d,}
and let they satisfy \eqref{RelationInit} and \eqref{eq_ini}. Then there exists a global in time weak solution to \eqref{S-appr} 
such that 
\eq{\label{est_kfixed}
&R_k^{-1},\ R_k,\ Q_k^{-1},\ Q_k,\ Z_k^{-1},\ Z_k\in L^\infty((0,T)\times\T^d),\\
& \pt R_k+\vu_k\cdot\Grad R_k\in L^\infty((0,T)\times\T^d),\\
& \pt Q_k+\vu_k\cdot\Grad Q_k\in L^\infty((0,T)\times\T^d),\\
 &\pt Z_k+\vu_k\cdot\Grad Z_k\in L^\infty((0,T)\times\T^d),\\
& \nabla_x \vu\in{L^\infty(0,T;BMO(\T^d))},\quad \Div\vu\in{L^\infty((0,T)\times\T^d)},\quad \partial_t \vu\in{L^\infty(0,T;L^2(\T^d))},
}
where Equations \eqref{S-appr1}, \eqref{S-appr2}, and \eqref{S-appr4} are satisfied in ${{\cal D}'((0,T)\times \T^d})$, and the initial conditions \eqref{ini_ep}   are satisfied in  ${\cal D}'(\T^d)$.
\end{thm}

The rest of this section is devoted to the proof of this Theorem. It will be divided into three main steps:
\begin{enumerate}
\item[Step 1.] Proof of existence of solutions to the Lagrangian reformulation of system \eqref{S-appr}.
\item[Step 2.] Construction of  characteristics for certain regularization of the flow.
\item[Step 3.] Passage to the limit with regularization parameter and retrieval of the Eulerian formulation.
\end{enumerate}

\subsection{Proof of Theorem \ref{approx} -- Step 1}

The starting point for this section is system \eqref{S-appr} written in the  Lagrangian coordinates. We omit the definition of these coordinates  at this level on purpose, we will come back to this issue in the consecutive steps of the proof of Theorem \ref{approx}.

We consider the following system for  $ (t,y)\in (0,T)\times\T^d$:
\begin{subequations}\label{S-L}
\begin{align}
 & \d_t r+ r \sigma
 =0,\label{S-L1}\\
 & \d_t q + q \sigma 
 =0, \label{S-L2}\\
  &{\sigma =  ({\cal T}_k(z))^{\gamma_+}-\{({\cal T}_k(z))^{\gamma_+}\}_{{\cal L}}},\label{S-L3}\\
 & q = \lr{1- \frac{r}{z}}z^\gamma,\quad r\leq z, \label{S-L4}
\end{align}
\end{subequations}
where by $\{\cdot\}_{{\cal L}}$ we denote the average on the torus 
\eq{\label{defL}
\{f\}_{{\cal L}} := \frac{1}{|\T^d|}\int_{\T^d} f(t,y) \exp \lr{\int_0^t \sigma (s,y) \, {\rm d}s}\dy.}
The unknowns of the \eqref{S-L} are $r=r(t,y)$, $q=q(t,y)$, $z=z(t,y)$,  and $\sigma=\sigma(t,y)$, whereas $z=z(r,q)$ is a unique solution to  \eqref{S-L4}, see the proof of Lemma \ref{Lemma_posit}.
We supplement the system \eqref{S-L} with the following initial conditions:
\eqh{
&r(0,y)=r_0(y),\quad q(0,y)=q_0(y),\quad 0< r_0,\ q_0<\infty\ \text{a.e. in } \T^d,\\
&r_0\leq z_0, \quad q_0 = \lr{1- \frac{r_0}{z_0}}z_0^\gamma.
}
{Note that there are} no space derivatives in system \eqref{S-L}, which is the main advantage of the Lagrangian reformulation. In fact it allows to transform the PDE system to the system of ODEs with a nonlocal term $\{({\cal T}_k(z))^{\gamma_+}\}_{{\cal L}}$. In addition to that, we know a-priori that  any solution of \eqref{S-L} satisfies the estimates
\begin{equation}\label{e1}
 \|\sigma\|_{L^\infty((0,T)\times\T^d)} \leq k^{\gamma_+},
\end{equation}
and 
\begin{equation}\label{e2}
\sup_{t\in(0,T)} \|r,\,q,\,r^{-1},\,q^{-1}\|_{L^\infty(\T^d)}(t) \leq  \|r_0,q_0,r_0^{-1},q_0^{-1}\|_{L^\infty(\T^d)} \exp(T k^{\gamma_+} ).
\end{equation}
This information suggests the choice of functional space for the fixed point theorem that will be used to prove the existence of solutions to system \eqref{S-L}. In what follows we show that the map 
\eqh{
\Phi\ :\  L^\infty((0,T)\times\T^d)\times L^\infty((0,T)\times\T^d)\to  L^\infty((0,T)\times\T^d)\times L^\infty((0,T)\times\T^d),}
\eqh{ \Phi(r,q)=(\bar r, \bar q),
}
where $\bar r, \bar q$ are the solutions to the following system in $(0,T)\times\T^d$
\begin{subequations}\label{Ban}
\begin{align}
 & \d_t \bar r+ r \sigma
 =0,\label{Ban1}\\
&  \d_t \bar q + q \sigma 
=0, \label{Ban2}\\
  &{\sigma =  ({\cal T}_k(z))^{\gamma_+}-\{({\cal T}_k(z))^{\gamma_+}\}_{\cal L}},\label{Ban4}\\
  &q = \lr{1-\frac{r}{z}} z^\gamma,\quad r\leq z, \label{Ban3}
\end{align}
\end{subequations}
is a contraction, at least for short time $t\in(0,T)$.
 
Let us denote 
\eqh{
 \delta \bar r = \bar r_1 - \bar r_2, \quad \delta \bar q = \bar q_1 - \bar q_2,
 \quad \delta r = r_1-r_2, \quad \delta q = q_1-q_2, \quad \delta \sigma = \sigma_1 -\sigma_2.
}
Using \eqref{Ban} we compute 
\eqh{
(\delta \overline r, \delta \overline q)=  \Phi(r_1,q_1)-\Phi(r_2,q_2),
}
we have
\begin{subequations}\label{Ban-dif}
\begin{align}
&  \d_t \delta \bar r+  \sigma_1 \delta r + r_2 \delta \sigma 
        =0, \label{Ban-dif1}\\
&  \d_t \delta \bar q + \sigma_1 \delta  q + q_2 \delta \sigma 
=0, \label{Ban-dif2} \\
 & \delta \sigma =  [({\cal T}_k(z_1))^{\gamma_+}-({\cal T}_k(z_2))^{\gamma_+}] -[\{({\cal T}_k(z_1))^{\gamma_+}\}_{\cal L} -\{({\cal T}_k(z_2))^{\gamma_+}\}_{\cal L}],\label{Ban-dif3}
\end{align}
\end{subequations}
where $z_1$ and $z_2$ are functions depending respectively on $(r_1,q_1)$ and $(r_2,q_2)$ through the nonlinear
implicit relation \eqref{Ban3}. Therefore, similarly to \eqref{pTR} we can show that
\eqh{
\|\partial_r {\cal T}_k(z)\|_{L^\infty((0,T)\times\T^d)}\leq C(k),\quad \mbox{and}\quad \|\partial_q {\cal T}_k(z)\|_{L^\infty((0,T)\times\T^d)}\leq C(k),
}
this implies, in particular that
\eq{\label{dsigmab}
\sup_{t\in(0,T)}\|\delta{\cal T}_k(z)\|_{L^\infty(\T^d)}(t)\leq C(k)\sup_{t\in(0,T)}\lr{\|\delta r\|_{L^\infty(\T^d)}(t) + \|\delta q\|_{L^\infty(\T^d)}(t)}.
}
To apply the Banach fixed point to the map $\Phi$, we need to show that  $\| \delta \bar r,\, \delta \bar q\|_{L^\infty((0,T)\times\T^d)}\leq C\| \delta  r,\, \delta  q\|_{L^\infty((0,T)\times\T^d)}$ with some constant $C<1$. To this purpose we need to estimate $\delta \sigma$ appearing in both equations, \eqref{Ban-dif1} and \eqref{Ban-dif2}, by $\delta r$ and $\delta q$. Analyzing equation \eqref{Ban-dif3}  we notice that since the first term on the r.h.s. can be treated using \eqref{dsigmab},  the only  challenge is due to the nonlocal term. Using \eqref{defL} and \eqref{dsigmab} we write
\eq{\label{e6}
&\left|\T^d| | \{({\cal T}_k(z_1))^{\gamma_+}\}_{\cal L}(t) -\{({\cal T}_k(z_2))^{\gamma_+}\}_{\cal L} (t)\right|\\
&=\left| \int_{\T^d} ({\cal T}_k(z(r_1,q_1)))^{\gamma_+} \exp\lr{\int_0^t \sigma_1(s,y) {\rm d}s}\dy \right.\\
&\hspace{5cm} \left.- \int_{\T^d} \lr{{\cal T}_k(z(r_2,q_2))}^{\gamma_+} \exp\lr{\int_0^t \sigma_2(s,y) {\rm d}s} \dy\right| \\
&\leq \left| \int_{\T^d} \exp\lr{\int_0^t \sigma_1(s,y) {\rm d}s} (({\cal T}_k(z(r_1,q_1)))^{\gamma_+}-({\cal T}_k(z(r_2,q_2)))^{\gamma_+}) \dy\right| \\
&\qquad- \left|\int_{\T^d} ({\cal T}_k(z(r_2,q_2)))^{\gamma_+} \left[ \exp\lr{\int_0^t \sigma_2(y,s) {\rm d}s}-  \exp\lr{\int_0^t \sigma_1(s,y) {\rm d}s}\right]\dy\right| \\
&\leq C(k)\lr{\|\delta r\|_{L^\infty(\T^d)}(t) + \|\delta q\|_{L^\infty(\T^d)}(t)}+ C(k) \int_{\T^d} \left|\int_0^t \delta \sigma (s,y) {\rm d}s\right|\dy \\
&\leq C(k)\lr{\|\delta r\|_{L^\infty(\T^d)}(t) + \|\delta q\|_{L^\infty(\T^d)}(t)} \\
&\qquad+ C(k)\,  t \sup_{s\in(0,t)} \left| \{({\cal T}_k(z_1))^{\gamma_+}\}_{\cal L}(s) -\{({\cal T}_k(z_2))^{\gamma_+}\}_{\cal L}(s)\right|,
}
where the last inequality follows from the last equation in \eqref{Ban-dif}. Taking now supremum over time $t\in (0,\tau)$ on both sides of \eqref{e6}  for $\tau C(k) \leq 1/2 |\T^d|$ we find
\begin{equation}\label{e7}
 \sup_{t\in(0, \tau)} | \{({\cal T}_k(z_1))^{\gamma_+}\}(t) -\{({\cal T}_k(z_2))^{\gamma_+}\}(t)|\leq
 C(k)\,  \sup_{t\in(0,\tau)} \lr{\|\delta r\|_{L^\infty(\T^d)}(t) + \|\delta q\|_{L^\infty(\T^d)}(t)}.
\end{equation}
With this estimate at hand we can return to (\ref{Ban-dif})  and compute  that 
\begin{equation}\label{e8}
 \sup_{0\leq t\leq \tau} \lr{\|\delta \bar r\|_{L^\infty(\T^d)}(t) + \|\delta \bar q\|_{L^\infty(\T^d)}(t)}\leq
 C(k)t \sup_{0\leq t\leq \tau}\lr{\|\delta r\|_{L^\infty(\T^d)}(t) + \|\delta q\|_{L^\infty(\T^d)}(t)}.
\end{equation}
So, choosing $\tau$ small such that   $\tau C(k) <1$,  map $\Phi$ is a contraction on a time interval $[0,\tau]$. Since $\tau$ depends only on the truncation parameter $k$ which is constant,  we can iterate this procedure to obtain a unique solution to  (\ref{S-L}) on the whole time interval $(0,T)$. 

\subsection{Proof of Theorem \ref{approx} -- Step 2}
\subsubsection{Explanation of the strategy}\label{Sect_strategy}
The existence of solutions proved before corresponds, in some sense, to the  Lagrangian reformulation of  system \eqref{S-appr}. We now want to define the Eulerian coordinates and show that we can recover 
the velocity vector field at the level of the Eulerian coordinates.  Using the mathematical jargon, we intend to solve the equation
\eq{\label{heur1}
\Div_x \vu(t,x)=\sigma(t,y),}
which means that for given $\sigma$, we will find $\vu$ and $x=x(t,y)$ such that the above equality is satisfied. Our candidate $x=x(t,y)$ is a solution to an ODE defining the Lagrangian transformation
\eq{\label{heur2}
\frac{d x}{dt} = \vu(t,x), \qquad x\vert_{t=0} = y.}
Combination of \eqref{heur1} and \eqref{heur2} leads to the nonlinear PDE-ODE system which we intend to solve using the following variant of the Leray-Schauder fixed point theorem.
\begin{thm}
Let $\Phi$ be a continuous, compact mapping, $X$ a Banach Space. Let for any $\zeta\in[0,1]$ the fixed point $v=\zeta\Phi(v)$, $v\in X$ be bounded. Then $\Phi$ possesses at least one fixed point in $X$.
\end{thm}

Let us now explain how we intend to apply the above theorem. We will consider a regularization of $\sigma$
\eqh{
 \sigma_\delta(t,y) = \sigma(t,y) \ast \kappa_\delta(y),
}
where  $\kappa_\delta$ is a standard mollifier and $\sigma$ is a known function -- the solution found in the previous step of the proof.  We further define a map $\Phi(\overline \vu) = \vu$
\eqh{
\Phi\ :\ {\cal C}(0,T; W^{2-\alpha,p}(\T^d))\to \ {\cal C}(0,T; W^{2-\alpha,p}(\T^d)),
}
for $\alpha>0$ arbitrary small, and {$p>1$}, in the following way: 

\medskip

\noindent  1. For a given $\overline \vu\in L^p(0,t; W^{2-\alpha,p}(\T^d))$ we use the Cauchy-Lipschitz theorem to find a unique $x=x(t,y)$ such that  
\eq{\label{ODE}
\frac{d x(t,y)}{dt} = \overline \vu(t,x(t,y)), \qquad x(t,y)\vert_{t=0} = y.}

\medskip

\noindent 2. We then differentiate \eqref{ODE} with respect to $y$, to check that $H(t,y) = \frac{\partial x}{\partial y}(t,y)$ satisfies the equation
$$\pt H(t,y)=\Grad_x \bar \vu(t,x(t,y)) H(t,y),\quad H(t,y)\vert_{t=0}=Id.$$
Therefore, integrating in time, we get:
\eqh{
\exp\lr{ - \int_0^T \|\nabla_x \bar \vu \|_{L^\infty(\T^d)} \dt}\leq \left\|H\right\|_{L^\infty((0,T)\times\T^d)} \leq \exp\lr{  \int_0^T \|\nabla_x \bar \vu \|_{L^\infty(\T^d)} \dt}.
}
Note also that the determinant of $H(t,y)$, $J(t,y)=det \,\frac{\partial x}{\partial y}(t,y)$ satisfies the equation
\eq{\label{Jf}
\pt J(t,y)=\Div_x\bar \vu(t,x(t,y)) J(t,y).
}
and so
\begin{equation}\label{j_est}
 \exp\lr{ - \int_0^T \|\Div_x \bar \vu\|_{L^\infty(\T^d)} \dt} \leq \|J\|_{L^\infty((0,T)\times\T^d)} \leq \exp\lr{ \int_0^T \|\Div_x \bar \vu\|_{L^\infty(\T^d)} \dt}.
\end{equation}
This means that $H(t,y)$ is invertible, moreover we have
\begin{equation}\label{f11}
 \left\|\frac{\partial y}{\partial x}\right\|_{L^\infty((0,T)\times\T^d)} \leq \exp\lr{ \int_0^T \|\nabla \bar \vu\|_{L^\infty(\T^d)} \dt}.
\end{equation}

\medskip

\noindent 3. Because $H$ is invertible we can express $y$ as a function of $t$ and $x$. For such $y=y(t,x)$ we will look for solution
\eq{\label{phid}
\vu(t,x)=\Grad\phi(t,x),\quad \mbox{where} \quad \lap_x \phi(t,x)=\sigma_\delta(t,y(t,x)).
}
\begin{rmk}
This approach guarantees not only that $\Div_x\vu(t,x)=\sigma_\delta(t,y(t,x))$, but also, that $\vu$ has a structure of a gradient flow and will satisfy \eqref{S-appr5} and \eqref{S-appr6}.
\end{rmk}
\begin{rmk}
The solution $\vu$ constructed above depends on the parameter $\delta$ and should be denoted $\vu^\delta$, but we omit this index at this stage of the proof.
\end{rmk}

\subsubsection{A priori estimates} 
We assume here that $\vu=\bar\vu \in {\cal C}(0,T; W^{2-\alpha,p}(\T^d))$ such that there exists $\vu(t,x)=\Grad_x \phi(t,x)$ and it satisfies
\eq{\label{phidl}
\lap_x\phi(t,x)= \zeta \sigma_\delta(t,y(t,x)),\quad {\zeta}\in[0,1].}
We will show that every solution of this equation is bounded in ${\cal C}(0,T; W^{2-\alpha,p}(\T^d))$ uniformly with respect to $\zeta $.

\bigskip

\noindent{\it Estimates of space derivative.} First note that the standard elliptic estimate for \eqref{phidl} gives
\eqh{\sup_{t\in(0,T)}\|\Grad_x\phi(t,x)\|_{W^{1,p}(\T^d)}\leq \zeta C(p)\|\sigma_\delta(t,y(t,x))\|_{L^\infty((0,T)\times \T^d)}
\leq C(p,k)
,}
and so 
\eq{\label{uBMO}
{\sup_{t\in(0,T)}\|\Grad\vu\|_{BMO(\T^d)}\leq C(k)}
}
for any $p<\infty$. In particular, taking $p>d$, we obtain $W^{1,p}(\T^d)\hookrightarrow \hookrightarrow C(\Ov{\T^d})$.
Differentiating equation \eqref{phidl} with respect to $x_k$ gives
\begin{equation}\label{f10}
 \Delta_x \frac{\partial\phi(t,x)}{\partial x_k} = \nabla_y \sigma_\delta(t,y(x,t)) \cdot \frac{\partial y}{\partial x_k},
\end{equation}
therefore using the estimate \eqref{f11} with $\vu=\bar\vu$ we obtain
\eqh{\sup_{t\in(0,T)} \|\nabla_x \vu\|_{W^{1,p}(\T^d)}
 \leq \zeta  C(p) \|\nabla_y \sigma_\delta\|_{L^\infty((0,T)\times\T^d)} \exp\lr{ T \sup_{t\in(0,T)} \|\nabla_x \vu\|_{L^\infty(\T^d)} } .}
We now use the following estimate for  $p >d$  
\begin{equation}\label{f12}
 \|\nabla f\|_{L^\infty(\T^d)} \leq C(p)(1+\|\nabla f\|_{BMO(\T^d)} \lr{1+\ln^+ (\|\nabla f\|_{W^{1,p}(\T^d)} + \|f\|_{L^\infty(\T^d)}) }^{1/2}).
\end{equation}
The proof of above inequality on can find in \cite{Ogawa} in Corollary 2.4 inequality (2.6). The original proof holds for the whole $\R^d$,  but taking an extension operator $E:\T^d \sim [0,1]^d \to \R^d$ preserving regularity in $W^{2,p}$ we get the same  result on the torus.
Together with \eqref{uBMO} and $ \zeta \leq 1$, we can write
\eqh{
 &\sup_{t\in(0,T)} \|\nabla_x \vu\|_{W^{1,p}(\T^d)} \\
 &
  \leq C(p, \delta) \exp\lr{ T \|\nabla_x \vu\|_{L^\infty(0,T; BMO(\T^d))}(1+\ln^+ \|\nabla_x \vu\|_{L^\infty(0,T; W^{1,p}(\T^d))})^{1/2}} \\
& \leq C(p,\delta) \exp\lr{ \frac{1}{2} \ln^+ \|\nabla_x \vu\|_{L^\infty(0,T; W^{1,p}(\T^d))} + C(1+T^2)} \\
& \leq C(p,\delta) e^{C(1+T^2)} \|\nabla_x \vu\|^{1/2}_{L^\infty(0,T; W^{1,p}(\T^d))}.
}
We therefore have
\begin{equation}\label{est_ud}
  \sup_{t\in(0,T)} \|\nabla_x\vu\|_{W^{1,p}(\T^d)}  \leq C.
\end{equation}
 In this way we find the {\it a-priori} information about $\vu$:
\eqh{
 \vu(t,x) \in L^\infty(0,T; W^{2,p}(\T^d)),\qquad \div_x \vu(t,x) = \sigma \ast \kappa_\delta (t,y(t,x)).}
The estimates uniform in $\delta$ are
\begin{equation}\label{est_u}
 \|\nabla_x \vu\|_{L^\infty(0,T;BMO(\T^d))} + \|\Div_x \vu\|_{L^\infty((0,T)\times\T^d)} \leq C.
\end{equation}

\bigskip

\noindent {\it Estimate  of the time derivative.} We want to check the time-regularity of $\vu(t,x)$. For fixed $\delta$ we expect better information, however here we present only the estimates uniform with respect to $\delta$. 
The therefore switch to the weak formulation  of (\ref{phidl}):
$$
\intO{ \Delta_x \phi (t,x) \pi(x)} = \zeta \intO{\sigma_\delta (t,y(x,t))\pi(x)}= 
\zeta \int_{\T^d}\sigma_\delta(t,y) \pi(x(t,y)) J(t,y) \, \dy,
$$
where $\pi$ is a smooth function on $\T^d$. We now differentiate this identity in time 
\eq{\label{ptsigma}
\Dt  \intO{\Delta_x \phi (t,x) \pi(x)} =&  
\zeta \int_{\T^d} \pt\left[\sigma_\delta(t,y) \pi(x(t,y)) J(t,y)\right]\, \dy\\
=& \zeta  \int_{\T^d}\pt\sigma_\delta(t,y)\pi(x(t,y)) J(t,y)\,\dy
+\zeta  \int_{\T^d}\sigma_\delta(t,y)\pt\pi(x(t,y)) J(t,y)\,\dy\\
&\quad+ \zeta  \int_{\T^d}\sigma_\delta(t,y)\pi(x(t,y)) \pt J(t,y)\,\dy.
}
We now need to estimate all terms on the r.h.s. of \eqref{ptsigma}. For the first one, we  differentiate \eqref{S-L4}
in time, use expressions \eqref{S-L1} and \eqref{S-L2} for $\pt r$ and $\pt q$, and proceed as in \eqref{dsigmab} to check  that
$$
|\pt \sigma_\delta |=|(\pt \sigma)\ast\kappa_\delta| \leq |\pt \sigma|  \leq C(k)( |r_t|+|q_t|) \in L^\infty((0,T)\times\T^d).
$$

Next, for the second term on the r.h.s. of \eqref{ptsigma} we use  \eqref{ODE} with $\bar\vu=\vu$ to write
$$
\pt \pi(x(t,y))= \nabla_x \pi \frac{dx}{dt}= \vu(t,x(t,y)) \nabla_x \pi(x(t,y)),
$$
so, thanks to \eqref{est_u}, $\pt \pi(x(t,y))$ is bounded in $L^\infty(0,T;L^2(\T^d))$.  At last, the formula \eqref{Jf} with $\bar\vu=\vu$ together with \eqref{j_est} provides that $\pt J(t,y)$ is bounded in $L^\infty((0,T)\times\T^d)$, so the third term on the r.h.s. of \eqref{ptsigma} is also bounded.

\medskip
Summarizing, we  get that $\lap_x\pt\phi$ is bounded in $L^\infty(0,T;W^{-1,2}(\T^d))$ uniformly with respect to $\delta$. In particular, using the Helmholtz decomposition, we deduce  that 
\eq{\label{est_ptu}
\|\pt\vu\|_{L^\infty(0,T;L^2(\T^d))}=\|\Grad_x\pt\phi\|_{ L^\infty(0,T;L^2(\T^d))}\leq C,} 
with the constant $C$ that does not depend on $\delta$.
 Obviously, combining the estimate \eqref{est_ptu} with \eqref{est_ud} and \eqref{est_u}, we verify that any fixed point satisfying \eqref{phidl} is uniformly bounded in ${\cal C}(0,T; W^{2-\alpha,p}(\T^d))$, independently of $\zeta $.

\subsubsection{The fixed point argument} 
We are now ready to proceed with the fixed point argument explained in Section \ref{Sect_strategy}.
We therefore take $\bar\vu \in {\cal C}(0,T;W^{2-\alpha,p}(\T^d))$, and show that the operator $\Phi(\bar\vu)=\vu$ defined through \eqref{phid} and \eqref{ODE} is continuous and compact in ${\cal C}(0,T;W^{2-\alpha,p}(\T^d))$.

Compactness is straightforward. Indeed, taking $\zeta =1$ in the system \eqref{phidl} we see that our a-priori estimates for $\vu$, estimates \eqref{est_ptu} with \eqref{est_ud} and \eqref{est_u}, stay in force. Hence, on account of the Aubin-Lions lemma, the map $\Phi$ is compact. 

We now check the continuity of the map $\Phi$, by investigating the difference of two solutions $\vu_1=\Grad_x\phi_1$ and $\vu_2=\Grad_x\phi_2$
\eq{\label{conv_phi}
\intO{ |\lap_x \phi_1(t,x) - \lap_x \phi_2(t,x)|^p} &= \intO{|\sigma_\delta(t,y_1(t,x)) - \sigma_\delta(t, y_2(t,x))|^p}\\
 &\leq C\|\Grad_y\sigma_\delta(t,y(t,x))\|^p_{L^\infty(\T^d)}\| y_1(t,x) -  y_2(t,x)\|^p_{L^\infty(\T^d)}.
}
Recalling \eqref{ODE}, we  note that for the regular characteristics intersecting in the point $(t,x)$ we have
\eqh{
 \|  y_1(x,t) -  y_2(x,t)\|_{L^\infty(\T^d)}\leq \int_0^t \|\bar\vu_1 - \bar\vu_2\|_{L^\infty(\T^d)} dt' .
}
So if $(\bar\vu_1-\bar\vu_2) \to 0$ in ${\cal C}(0,T;W^{2-\alpha,p}(\T^d))$, then of course
$(\bar\vu_1-\bar\vu_2) \to 0$ in ${\cal C}(0,T;L^{\infty}(\T^d))$, which thanks to \eqref{conv_phi} implies that 
$(\phi_1 - \phi_2) \to 0$ in ${\cal C}(0,T;W^{2,p}(\T^d))$. Therefore, using boundedness of $\vu_1,\ \vu_2$ in  $L^\infty(0,T;W^{2,p}(\T^d))$, and interpolation, we deduce that $(\vu_1-\vu_2) \to 0$ in $L^p(0,T;W^{2-\alpha,p}(\T^d))$. By Leray-Schauder fixed point theorem we have the existence for approximative system for $\delta$ fixed.

\subsection{Proof of Theorem \ref{approx} -- Step 3}
 We want to let  $\delta\to 0$ in the equation 
\begin{equation}\label{f20}
 \Delta_x \phi^\delta(t,x^\delta) = \sigma_\delta(t,y(t,x^\delta))
\end{equation}
remembering that  $\vu^\delta(t,x^\delta)=\Grad_x\phi^\delta(t,x^\delta)$ and $x^\delta $ is associated with the flow $\vu^\delta$ via \eqref{ODE} with $\bar\vu=\vu^\delta$. We know that uniformly with respect to $\delta$ we have
\eqh{
 \|\nabla_x \vu^\delta\|_{L^\infty(0,T;BMO(\T^d))} + \|\Div_x \vu^\delta\|_{L^\infty((0,T)\times\T^d)}  + \|\partial_t \vu^\delta\|_{L^\infty(0,T;L^2(\T^d))} \leq C. 
}
We use again the weak form of (\ref{f20}), we have
\eq{\label{weak_sold}
 \intT{\!\!\!\int_{\T^d}\Delta_x \phi^\delta (t,x^\delta) \xi (t,x^\delta) \dx^\delta} = \intT{\!\!\!\int_{\T^d} \sigma \ast \kappa_\delta (t,y(x^\delta,t)) \xi(t,x^\delta)\dx^\delta} \\
 = \intT{\!\!\!\int_{\T^d} \sigma \ast \kappa_\delta (t,y) \xi(t,x^\delta(t,y)) J^\delta(t,y) \dy},
}
for any smooth $\xi$, where we denoted
$$J^\delta(t,y)=\exp\lr{ \int_0^t \Div_x \vu^\delta(t',x^\delta(t',y)) {\rm d}t'}=\exp\lr{ \int_0^t \sigma \ast \kappa_\delta (t',y) {\rm d}t'}.$$
We see $\phi^\delta \to \phi$ weakly in $L^\infty(0,T;W^{2,p}(\T^d))$, and
$\sigma\ast \kappa_\delta \to \sigma$ a.e. in $(0,T)\times\T^d$.

Using Crippa-Dellelis result  from \cite[Theorem 2.9 (stability of the flow)]{CrDe} saying
\eqh{
 \sup_{t\in(0,T)} \| x(t,y) -  x^\delta(t,y)\|_{L^1(\T^d)} \leq C | \ln ( \| \vu -\vu^\delta\|_{L^1((0,T) \times \T^d)})|^{-1}
}
we therefore get
$$
x^\delta(t,y) \to x(t,y) \mbox{ in } L^\infty(0,T;L^1(\T^d)),
$$
since $\vu^\delta$ is compact in $L^1((0,T)\times \T^d)$.
Hence we can let $\delta\to0$ in both sides of \eqref{weak_sold} to obtain
\eqh{
 \intTO{\Delta_x \phi (t,x) \xi (t,x)} =  \intT{\!\!\!\int_{\T^d} \sigma  (t,y) {\xi(t,x(y,t)) }J(t,y)\, \dy.}
}
Now we substitute for $\sigma$ using \eqref{Ban4}, and use Lemma 3.1. from \cite{CoCrSp} to write the weak formulation of the momentum equation in the Eulerian coordinates
{\eq{&\intTO{\Div_x \vu (t,x) \xi (t,x)}\label{relation} \\
&\quad= \intTO{\Big(({\cal T}_k(Z(t,x)))^{\gamma_+}-\left\{({\cal T}_k(Z(t,x)))^{\gamma_+}\right\}\Big)\xi(t,x)},}}
where $Z(t,x)$ is defined by $Z(t,x(t,y))=z(t,y)$. Defining $R(t,x(t,y))=r(t,y)$ and $Q(t,x(t,y))=q(t,y)$ we can also pass to the limit in the weak formulation of \eqref{S-appr1} and \eqref{S-appr2}. Note that for the limit case, characteristics $x(t,y)$ are now  well defined. So, solutions given in the Lagrangian coordinates give us weak solutions to the original approximate problem (\ref{S-appr}). 

\section{Existence of solutions to system \eqref{S}}
Having proven the existence of solutions to the approximate system we now intend to pass to the limit with the approximation parameter $k\to\infty$,  recovering the system \eqref{S} and concluding the proof of Theorem \ref{Thm:main}. This section will also include rigorous justification of the compactness result from the previous section at the level of the approximate system.

From Theorem \ref{approx} it follows that there exists a sequence $\{R_k,Q_k,Z_k,\vu_k\}_{k=1}^\infty$ satisfying  system \eqref{S-appr} in the sense of distributions belonging to a class \eqref{est_kfixed}.
Note, however, that this is not information uniform with respect to $k$. 
\subsection{Uniform estimates}
Although \eqref{est_kfixed} is not uniform with respect to $k$, we can still use \eqref{e2} from the previous section to deduce that
$$R_k\geq 0,\quad Q_k\geq 0,$$
uniformly with respect to $k$. Moreover, repeating the proof  of Lemma \ref{Lemma_posit}, we find unique $Z_k$ defined via \eqref{S-appr3}, therefore uniformly w.r.t. $k$ we also have that
$$R_k\leq Z_k.$$
In the rest of this section we will skip the index $k$ where no confusion can arise.

 \subsubsection{Estimates of $({\cal T}_k(Z))^{\gamma^+}$ uniform w.r.t.  $k$}
 
We use the renormalized equation for $Z$, \eqref{SZ}, that can be derived from \eqref{S-appr1} and \eqref{S-appr2} using \eqref{S-appr3}. Testing this equation respectively  by $(\gamma_+-1) ({\cal T}_k(Z))^{\gamma_+-2} {\cal T}'_k(Z)$ and by $(\gamma_+-\gamma) ({\cal T}_k(Z))^{\gamma_+-\gamma-1} {\cal T}'_k(Z)$, we obtain
\eqh{
&\partial_t ({\cal T}_k(Z))^{\gamma^+-1} + \vu \cdot \nabla ({\cal T}_k(Z))^{\gamma^+-1} 
     + (\gamma^+ -1) ({\cal T}_k(Z))^{\gamma^+-2}{\cal T}'_k(Z) Z\Div \vu\\
&  \hskip3cm    +  (\gamma^+ -1) \frac{(1-\gamma)(Z-R) ({\cal T}_k(Z))^{\gamma^+-2}{\cal T}'_k(Z)Z  }{\gamma (Z-R) + R} \Div \vu
= 0,
}
and
\eqh{
&\partial_t ({\cal T}_k(Z))^{\gamma^+-\gamma} + \vu \cdot \nabla ({\cal T}_k(Z))^{\gamma^+-\gamma} 
     + (\gamma^+ -\gamma) ({\cal T}_k(Z))^{\gamma^+-\gamma-1}{\cal T}'_k(Z)Z \Div \vu\\
&  \hskip3cm    +  (\gamma^+ -\gamma) \frac{(1-\gamma)(Z-R) ({\cal T}_k(Z))^{\gamma^+-\gamma-1}{\cal T}'_k(Z)Z  }{\gamma (Z-R) + R} \Div \vu 
= 0.
}
After simplification we get
\eqh{
&\partial_t ({\cal T}_k(Z))^{\gamma^+-1} + \vu \cdot \nabla ({\cal T}_k(Z))^{\gamma^+-1} 
     \\
&  \hskip3cm    +  (\gamma^+ -1) \frac{Z^{\gamma^+}  }{\gamma (Z-R) + R} \Div \vu 1_{Z\leq k}
= 0,
}
and
\eqh{
&\partial_t ({\cal T}_k(Z))^{\gamma^+-\gamma} + \vu \cdot \nabla ({\cal T}_k(Z))^{\gamma^+-\gamma} \\
&  \hskip3cm    +  (\gamma^+ -\gamma) \frac{Z^{\gamma^+-\gamma+1}}{\gamma (Z-R) + R} \Div \vu1_{Z\leq k}
= 0.
}
If we multiply the first equation by $R$ and the second by $Q$ and use the continuity equations \eqref{S-appr1} and \eqref{S-appr2} we get
\eq{\label{RTk}
&\partial_t \lr{R ({\cal T}_k(Z))^{\gamma^+-1}} + \Div\lr{R ({\cal T}_k(Z))^{\gamma^+-1} \vu}
     \\
&  \hskip3cm    +  (\gamma^+ -1) \frac{RZ^{\gamma^+}  }{\gamma (Z-R) + R} \Div \vu 1_{Z\leq k}
       =0,
}
and
\eq{\label{QTk}
&\partial_t \lr{Q ({\cal T}_k(Z))^{\gamma^+-\gamma}} + \Div(Q \lr{{\cal T}_k(Z))^{\gamma^+-\gamma} \vu} \\
&  \hskip3cm    +  (\gamma^+ -\gamma) \frac{QZ^{\gamma^+-\gamma+1}}{\gamma (Z-R) + R} \Div \vu 1_{Z\leq k}
       =0.
}
Multiplying equations \eqref{RTk}, \eqref{QTk} by $\frac{1}{(\gamma^+-1)}$, $\frac{1}{(\gamma^--1)}$, respectively and noticing that $\frac{\gamma^+-\gamma}{(\gamma^--1)}=\gamma$ we obtain 
\eqh{
&\frac{1}{(\gamma^+-1)} \pt\lr{ R \, ({\cal T}_k(Z))^{\gamma^+-1}}  + \frac{1}{(\gamma^--1)} \pt\lr{  Q ({\cal T}_k(Z))^{\gamma^+-\gamma} }\\
&+\frac{1}{(\gamma^+-1)} \Div\lr{R ({\cal T}_k(Z))^{\gamma^+-1} \vu}+\frac{1}{(\gamma^--1)} \Div\lr{Q( {\cal T}_k(Z))^{\gamma^+-\gamma} \vu} \\
&   + Z^{\gamma^+}\frac{R+\gamma Q Z^{1-\gamma}}{\gamma (Z-R) + R}\Div \vu 1_{Z\leq k}
\leq 0,
}
which, after integrating with respect to space, and noting that $Q Z^{1-\gamma} = Z- R$, gives
\eqh{
&\frac{1}{(\gamma^+-1)} \Dt\intO{ R \, ({\cal T}_k(Z))^{\gamma^+-1}}  + \frac{1}{(\gamma^--1)} \Dt \intO{  Q ({\cal T}_k(Z))^{\gamma^+-\gamma} }\\
   &\quad+ \intO{ \Div \vu\,  Z^{\gamma^+}1_{Z\leq k}}
   \leq 0.
}
Let us note that $Z^{\gamma^+}1_{Z\leq k}=({\cal T}_k(Z))^{\gamma^+}$ and $\left\{{\cal T}_k(Z))^{\gamma^+}\right\}$ is constant in space, therefore
$$\intO{ \Div \vu\,  Z^{\gamma^+}1_{Z\leq k}}=\intO{ \Div \vu\,  \lr{Z^{\gamma^+}1_{Z\leq k}-\left\{{\cal T}_k(Z))^{\gamma^+}\right\}}}.$$
And so, using density of smooth functions in $L^2$, and taking in  \eqref{relation} $\xi=\Div\vu$, we get
\eq{\label{en_approx}
&\frac{1}{(\gamma^+-1)} \Dt\intO{ R \, ({\cal T}_k(Z))^{\gamma^+-1}}  + \frac{1}{(\gamma^--1)} \Dt \intO{  Q ({\cal T}_k(Z))^{\gamma^+-\gamma} }\\
   &\quad + \intO{\Div^2 \vu} \leq 0,
}
hence, 
\begin{equation}\label{en_approx-2}
\Dt \left[ \int { \frac{1}{(\gamma^+-1)} { R \, ({\cal T}_k(Z))^{\gamma^+-1}}  + \frac{1}{(\gamma^--1)} {  Q ({\cal T}_k(Z))^{\gamma^+-\gamma} }}\right]\leq 0.
\end{equation}
We now observe that
\eq{
&R  ({\cal T}_k(Z))^{\gamma^+-1}+Q ({\cal T}_k(Z))^{\gamma^+-\gamma}\\
&=R  ({\cal T}_k(Z))^{\gamma^+-1}+(Z-R)Z^{\gamma-1} ({\cal T}_k(Z))^{\gamma^+-\gamma}\\\
&= Z^{\gamma} ({\cal T}_k(Z))^{\gamma^+-\gamma} + R[k^{\gamma-1}-Z^{\gamma-1}]1_{Z\geq k} ({\cal T}_k(Z))^{\gamma^+-\gamma}.
}
Since $\gamma-1<0$, the last term is nonnegative, therefore 
\eq{\label{ineq}
({\cal T}_k(Z))^{\gamma^+} \le 
 R  ({\cal T}_k(Z))^{\gamma^+-1}+Q ({\cal T}_k(Z))^{\gamma^+-\gamma}.
 }
Therefore, using \eqref{ineq} and integrating \eqref{en_approx-2} over time, we obtain
\eq{\label{Tk!!!}
\|T_k(Z)\|_{L^\infty(0,T;L^{\gamma^+}(\T^d))}\leq C,
}
uniformly w.r.t. $k$. 

\subsubsection{Uniform estimates of $R^{\gamma^+}$ and on $Z^{\gamma^+}$}
Recall that the equation on $R$ reads
$$\partial_t R + \Div (R\vu) 
     = 0,$$
and is satisfied in the sense of distributions on $(0,T)\times\T^d$. Due to \eqref{est_kfixed}, this equation is satisfied in the renormalized sense, therefore we have
$$\partial_t R^{\gamma^+} + \Div (R^{\gamma^+} \vu) 
       + (\gamma^+-1)\, R^{\gamma^+} \Div\vu 
     = 0.$$
{ Using \eqref{S-appr4} and integrating} with respect to space, it gives
 \eqh{
 & \frac{d}{dt} \intO{R^{\gamma^+}} + (\gamma^+-1) \intO{R^{\gamma^+} ({\cal T}_k(Z))^{\gamma^+}}\\
& \hskip3cm =
      (\gamma^+-1) \intO{R^{\gamma^+}} \intO{({\cal T}_k(Z))^{\gamma^+}}. 
 } 
 Using \eqref{Tk!!!}, we get  through Gronwall Lemma 
\eq{\label{R!!!}
\|R\|_{L^\infty(0,T;L^{\gamma^+}(\T^d))}+\|R^{\gamma^+} ({\cal T}_k(Z))^{\gamma^+}\|_{L^1((0,T)\times\T^d)}\leq C,
}
 uniformly w.r.t. $k$. Performing exactly the same procedure for $Z$ we obtain
 \eq{\label{Z!!!}
\|Z\|_{L^\infty(0,T;L^{\gamma^+}(\T^d))}+\|{\cal T}_k(Z)\|_{L^{2\gamma^+}((0,T)\times\T^d)}\leq C.
}
 These bounds are important to pass to the limit in the system (\ref{S-appr}) with respect
 to $k$. Note that the bound on $Z$ provides the same bound on
 $({\cal T}_k(Z))$. 
Finally, renormalizing the approximate equation for $Q$ we easily deduce that
\eq{\label{Q!!!}
\|Q\|_{L^\infty(0,T;L^{\gamma^-}(\T^d))}+\|Q^{\gamma^-} ({\cal T}_k(Z))^{\gamma^+}\|_{L^1((0,T)\times\T^d)}\leq C.
}
\subsubsection{Uniform estimates of $\vu$}
With the estimate \eqref{Z!!!} at hand we can now estimate $\Div\vu$ directly from \eqref{S-appr4} , we have 
$$\|\Div\vu\|_{L^\infty(0,T;L^{1}(\T^d))}+\|\Div\vu\|_{L^{2}((0,T)\times\T^d)}\leq C.$$
Using the fact that ${\rm rot}\vu=0$, we therefore find that
 \eq{\label{v!!!}
\|\vu\|_{L^{2}(0,T; W^{1,2}(\T^d))}\leq C.}

\subsection{Compactness argument}
The uniform estimates from the previous section are sufficient to perform the limit passage $k\to\infty$ in all linear terms of the approximate system \eqref{S-appr}. Estimating the time derivatives of $R_k$ and $Q_k$ from equations \eqref{S-appr1}, \eqref{S-appr2}, respectively,
and using the uniform estimate \eqref{v!!!} together with a Div-Curl type argument, we justify that $R_k\vu_k$ converges to $R\vu$ and $Q_k\vu_k$ converges to $Q\vu$ in the sense of distributions.

The last problem to solve is to pass to the limit  in the nonlinear term in the momentum equation \eqref{S-appr4} -- the pressure. This requires a strong convergence of the sequence $\{Z_k\}_{k=1}^\infty$ approximating $Z$. For the moment we only know that $Z_k$ converges weakly$^*$ to $Z$ in ${L^\infty(0,T;L^{\gamma^+}(\T^d))}$. To improve this convergence, we adapt the compactness criterion presented in Section \ref{Sect:3}. 

We first justify how to obtain the equivalent of equations \eqref{eqR} and \eqref{eqZ} on the approximate level. First, we write the equation for $R_x-R_y$ using \eqref{S-appr1} we get
\eqh{
&\pt (R_x-R_y) + \Div_x (\vu_x\lr{R_x-R_y}) +
\Div_y (\vu_y\lr{R_x-R_y})  \\
&=\frac 12 (\Div_x \vu_x + \Div_y \vu_y) \lr{R_x-R_y} \\
&\quad - \frac 12 (\Div_x \vu_x-\Div_y \vu_y)(R_x+R_y).
}
Regularizing this equation over the space variables, i.e. testing the equation by $\xi_\eta(x-\cdot)\xi_\eta(y-\cdot)$, where $\xi_\eta$ is a standard mollifier, and denoting $S_\eta[f]=\xi\ast_y(\xi\ast_x f)$ we obtain
\eq{\label{eqRdiff}
&\pt S_\eta[R_x-R_y] + \Div_x (\vu_xS_\eta[R_x-R_y]) +
\Div_y (\vu_yS_\eta[R_x-R_y])  \\
&=r_\eta+\frac 12 S_\eta\left[(\Div_x \vu_x + \Div_y \vu_y) \lr{R_x-R_y} \right]\\
&\quad - \frac 12 S_\eta\left[(\Div_x \vu_x-\Div_y \vu_y)(R_x+R_y)\right],
}
where on account of the Friedrichs commutator lemma, $r_\eta\to 0$ in $L^p((0,T)\times\T^{2d})$ for any $p<\infty$.
We then multiply \eqref{eqRdiff} by
\eqh{
S_\eta[R_x-R_y] \lr{S_\eta[R_x-R_y] ^2+\delta}^{\frac{\beta-2}{2}},
}
where $\delta>0$. Note that on account of \eqref{est_kfixed} this function is bounded in $L^\infty((0,T)\times\T^d)$ uniformly w.r.t. $\eta$. 
We therefore obtain
\eq{\label{eqR_ren}
&\frac{1}{\beta}\pt\lr{ \lr{S_\eta[R_x-R_y]}^2+\delta}^{\frac\beta2}\\
&\quad +\frac{1}{\beta} \Div_x \lr{\lr{ \lr{S_\eta[R_x-R_y]}^2+\delta}^{\frac\beta2}\vu_x} +
\frac{1}{\beta}\Div_y \lr{\lr{ \lr{S_\eta[R_x-R_y]}^2+\delta}^{\frac\beta2}\vu_y}  \\
&=\Big\{r_\eta+\frac 12 S_\eta\left[(\Div_x \vu_x + \Div_y \vu_y) \lr{R_x-R_y} \right]- \frac 12 S_\eta\left[(\Div_x \vu_x-\Div_y \vu_y)(R_x+R_y)\right]\Big\}\\
       &\qquad\times S_\eta[R_x-R_y] \lr{S_\eta[R_x-R_y] ^2+\delta}^{\frac{\beta-2}{2}}.
}
This equation can now be multiplied by $K_h(x-y)\lr{w_x+w_y}$, where $K_h$ is the same as in Section \ref{ssec:prelim}, and $w_x=w(t,x)$ satisfies  the equation
\eq{\label{def_w_ep}
\left\{\begin{array}{l}
\pt w+\vu\cdot\Grad w+ \theta {\cal D}_k w=0,\\
w(0,x)=1,
\end{array}\right.
}
with
\eq{\label{def_D_ep}
{\cal D}_k=M|\Grad\vu|+|\Div\vu|+({\cal T}_k(Z))^{\gamma^+}+ \{({\cal T}_k(Z))^{\gamma^+}\}.
}
The choice of ${\cal D}_k$ is to accommodate the extra terms appearing in the main compactness estimate, and will be explained later,
the same applies to the choice of the constant $\theta$. The most important observation at this level is that for any $k$, 
fixed ${\cal D}_k\in L^\infty((0,T)\times\T^d)$ and that due to regularity of   $\vu=\vu_k$, see \eqref{est_kfixed}, Equation \eqref{def_w_ep} 
has a unique distributional solution, see for example Lemma 6.10 and Lemma 6.12 from \cite{BJ}.

Therefore, letting  $\delta\to 0$, $\beta\to 1$, and $\eta\to 0$ in {$\int_{\T^{2d}}$\eqref{eqR_ren} $\dx$} and using the dominated convergence theorem, we recover \eq{\label{eqRdiff}
&\Dt\iintO{ K_h(x-y)|R_x-R_y|\lr{w_x+w_y}}  \\
&=\iintO{ \Grad K_h(x-y)(\vu_x-\vu_y)|R_x-R_y| (w_x+w_y)}\\
&\quad - \iintO{K_h(x-y)(\Div_x \vu_x-\Div_y \vu_y)R_xs_Rw_x}\\
&\quad+2 \iintO{K_h(x-y)|R_x-R_y| \lr{\pt w_x+\vu_x\cdot\Grad w_x+\Div_x\vu_x w_x}}.
}
To obtain the equation for $|Z_x-Z_y|$ we proceed in a similar manner. As a consequence, the analogue of \eqref{Rt2} we can be now written as
 \eq{\label{Sep}
&\Dt S(t)= \iintO{ \Grad K_h(x-y)(\vu_x-\vu_y)O_{x-y}  (w_x+w_y)}\\
&= \iintO{ \Grad K_h(x-y)(\vu_x-\vu_y)O_{x-y} (w_x+w_y)}\\
&\quad - \iintO{K_h(x-y)(\Div_x \vu_x-\Div_y \vu_y)[R_xs_R+Z_x s_Z]w_x}\\
&\quad+2(1-\gamma)\iintO{K_h(x-y)\left[\frac{ (Z_y-R_y)Z_y\Div_y\vu_y}{ \gamma(Z_y-R_y)+R_y}-
\frac{ (Z_x-R_x)Z_x\Div_x \vu_x}{ \gamma(Z_x-R_x)+R_x} 
\right] s_Zw_x}\\
&\quad+2 \iintO{K_h(x-y)O_{x-y} \lr{\pt w_x+\vu_x\cdot\Grad w_x+\Div_x\vu_x w_x}}\\
&=\sum_{i=1}^4 I_i,
}
where $O_{x-y}=|R_x-R_y|+ |Z_x-Z_y|$.
For abbreviation we will only comment on the changes due to the presence of new approximation terms and truncations. 
   The changes due to truncation in the momentum equation apply to terms $I_2$, and $I_3$. For the first one we use {\eqref{S-appr4}  to write}
\eq{
 &-\lr{ \Div_x\vu_x-\Div_y\vu_y}[R_xs_R+Z_x s_Z]\\
& =- \lr{ ({\cal T}_k(Z_x))^{\gamma^+}-({\cal T}_k(Z_y))^{\gamma^+}}R_xs_R
 - \lr{ ({\cal T}_k(Z_x))^{\gamma^+}-({\cal T}_k(Z_y))^{\gamma^+}}Z_x s_Z,
 }
 that gives the following contribution to the l.h.s. of \eqref{Sep} 
\eq{\iintO{K_h(x-y)\left|({\cal T}_k(Z_x))^{\gamma^+}-({\cal T}_k(Z_y))^{\gamma^+}\right|(-\alpha_x)Z_xw_x}.
}
For $I_3$ we write as in \eqref{re1}
\eq{\label{re1ap}
&(1-\gamma)\left[\frac{ (Z_y-R_y)Z_y\Div_y\vu_y}{ \gamma(Z_y-R_y)+R_y}-
\frac{ (Z_x-R_x)Z_x\Div_x \vu_x}{ \gamma(Z_x-R_x)+R_x} 
\right]s_Z\\
&=(1-\gamma)\left[\frac{ (Z_y-R_y)Z_y({\cal T}_k(Z_y))^{\gamma^+}}{ \gamma(Z_y-R_y)+R_y}-
\frac{ (Z_x-R_x)Z_x({\cal T}_k(Z_x))^{\gamma^+}}{ \gamma(Z_x-R_x)+R_x} 
\right]s_Z\\
&(1-\gamma)\left[ \frac{ (Z_x-R_x)Z_x}{ \gamma(Z_x-R_x)+R_x}-
\frac{(Z_y-R_y)Z_y}{ \gamma(Z_y-R_y)+R_y} 
\right] s_Z \{({\cal T}_k(Z))^{\gamma^+}\}\\
&=(1-\gamma)\frac{(Z_y-R_y)}{ \gamma(Z_y-R_y)+R_y}\lr{Z_y({\cal T}_k(Z_y))^{\gamma^+}-Z_x({\cal T}_k(Z_x))^{\gamma^+}} s_Z\\
&+(\gamma-1)
\left[
\frac{ (Z_x-R_x)}{ \gamma(Z_x-R_x)+R_x}-
\frac{(Z_y-R_y)}{ \gamma(Z_y-R_y)+R_y}
\right] Z_x({\cal T}_k(Z_x))^{\gamma^+} s_Z\\
&-(\gamma-1)\frac{(Z_y-R_y)}{ \gamma(Z_y-R_y)+R_y}\lr{Z_x-Z_y} s_Z \{({\cal T}_k(Z_x))^{\gamma^+}\} \\
&-(\gamma-1)
\left[
\frac{ (Z_x-R_x)}{ \gamma(Z_x-R_x)+R_x}-
\frac{(Z_y-R_y)}{ \gamma(Z_y-R_y)+R_y}
\right] Z_x\{({\cal T}_k(Z_x))^{\gamma^+} \}s_Z.
}
The first term has a good sign, while the rest of terms can be treated exactly as in the Section \ref{Sect:3}, using in particular the  definition of the weight \eqref{def_w_ep} and \eqref{def_D_ep} with $\theta$ sufficiently large  to absorb all the unwanted terms  by  $I_4$ in \eqref{Sep}.

The above considerations allow us to deduce 
\eqh{\label{without_w_ep}
\limsup_k\lr{\frac{1}{|\log h_0|}\iintO{ \calK_{h_0}(x-y)\Big[ \left|(Z_x)_k-(Z_y)_k\right|+ |(R_x)_k-(R_y)_k|\Big] }}\to 0,
}
by noticing that the thesis of Proposition \ref{prop:w} stay in force also for $(X,\vu)$ satisfying the continuity equation with a friction term.
Indeed, since ${\cal D}_k\in L^2((0,T)\times\Omega)$ independently of $k$ (see \eqref{Z!!!} and \eqref{v!!!}),  the crucial estimate \eqref{boundlogw} 
holds for both sequences $\{R_k\}_{k=1}^\infty$, and $\{Z_k\}_{k=1}^\infty$ independently of $k$. This observation finishes the proof of Theorem \ref{Thm:main}. $\Box$
  With this at hand, using Remark \ref{Rmk_1} to define $\alpha^+$ and $\alpha^-$, we find a weak solution to our original system \eqref{spec}. 
The proof of Theorem \ref{Mainbifluid} is therefore complete. $\Box$

\bigskip 

{\bf Acknowledgments.}  The author D.B. is partly supported by the 
ANR- 13-BS01- 0003-01 project DYFICOLTI, by the ANR-16-CE06-0011-02 
FRAISE and by the project  TelluS-INSMI-MI (INSU) CNRS. The  authors P.B.M. and E.Z. have been partly supported by National Science Centre grant 2014/14/M/ST1/00108 (Harmonia).
{The authors want to thank P.E. Jabin and the anonymous  referee for their valuable remarks about the paper.}

\end{document}